\documentclass[12pt]{amsart}
\usepackage{a4wide,amsfonts,amssymb,stmaryrd,amscd,amsmath,latexsym,amsbsy}
\numberwithin{equation}{section}
\theoremstyle{plain}
\newtheorem{thm}{Theorem}[section]

\newtheorem{prop}[thm]{Proposition}

\newtheorem{lem}[thm]{Lemma}
\newtheorem{cor}[thm]{Corollary}

\theoremstyle{remark}
\newtheorem{rema}[thm]{Remark}
\title{Hyperbolic beta integrals} 
\author{Jasper V. Stokman}
\address{J.V. Stokman,
KdV Institute for Mathematics, Universiteit van Amsterdam,
Plantage Muidergracht 24, 1018 TV Amsterdam, The Netherlands.}
\email{jstokman@science.uva.nl}
\begin{document}
\begin{abstract}
Hyperbolic beta integrals are analogues of Euler's beta integral in which
the role of Euler's gamma function is taken over by Ruijsenaars' 
hyperbolic gamma function. 
They may be viewed as $(q,\widetilde{q}\,)$-bibasic analogues
of the beta integral in which the two bases $q$ and $\widetilde{q}$ are
interrelated by modular inversion, and they entail 
$q$-analogues of the beta integral for $|q|=1$.
The integrals under consideration are the hyperbolic analogues of the
Ramanujan integral, the Askey-Wilson integral and
the Nassrallah-Rahman integral. We show that the hyperbolic
Nassrallah-Rahman integral is a formal limit case 
of Spiridonov's elliptic Nassrallah-Rahman integral.
\end{abstract}
\maketitle

\section{Introduction} 
\noindent
Euler's gamma function is defined by 
\[\Gamma(z)=\int_0^\infty x^{z-1}e^{-x}dx,\qquad \hbox{Re}(z)>0.
\]
In the fundamental paper \cite{R1}, Ruijsenaars defined 
gamma functions of rational, trigonometric, 
hyperbolic and elliptic type, Euler's gamma function being of rational
type. Accordingly, one expects to have extensive theories
on special functions of rational,
trigonometric, hyperbolic and elliptic type. The rational and trigonometric
special functions are precisely the special functions of hypergeometric and
basic hypergeometric type, which have been thoroughly
studied (see e.g. \cite{AAR} and \cite{GR}). Systematic studies of 
hyperbolic and elliptic special functions have
commenced only recently, see e.g. \cite{NU}, \cite{R2}, \cite{S} for the 
hyperbolic case and \cite{FT}, \cite{Sp}, \cite{SZ} for the elliptic case.

A basic step in the development of special functions of a given type
is the derivation of the associated beta integrals. 
The beta integral of rational type is Euler's beta integral
\begin{equation}\label{betaEuler}
\int_0^1 x^{a-1}(1-x)^{b-1}dx=\frac{\Gamma(a)\Gamma(b)}{\Gamma(a+b)},
\qquad \hbox{Re}(a), \hbox{Re}(b)>0.
\end{equation}
Trigonometric and elliptic analogues of beta integrals (involving
the trigonometric and the elliptic gamma function respectively)
have been studied in detail. The goal of this paper is to derive 
hyperbolic analogues of beta integrals. 

The importance of beta integrals lies in its variety of applications.
The rational beta integral is the normalization constant for the orthogonality
measure of the Jacobi polynomials and, more generally, multivariate
analogues of the beta integral arise as normalization constants in the 
theory of zonal spherical functions on compact symmetric spaces.
Beta type integrals also appear in number theory, combinatorics, 
conformal field theories and in certain completely integrable systems. 

Trigonometric beta integrals have been studied extensively,
see e.g. \cite{A2}, \cite{A3}, \cite{AR}, \cite{AW}, \cite{NR}.
In this case the role of Euler's gamma function, 
as well as of functions like
$x^{a-1}$ and $(1-x)^{b-1}$, are taken over by (quotients of)
Ruijsenaars' \cite{R1} trigonometric gamma functions, or equivalently, 
(quotients of) $q$-Pochhammer symbols. 
The $q$-Pochhammer symbol is defined for $|q|<1$ by
\[\bigl(a;q\bigr)_{\infty}=\prod_{j=0}^{\infty}(1-aq^j),
\]
and products of $q$-Pochhammer symbols will be denoted by
\begin{equation}\label{productPoch}
\begin{split}
\bigl(a_1,\ldots,a_m;q\bigr)_{\infty}&=\prod_{j=1}^m\bigl(a_j;
q\bigr)_{\infty},\\
\bigl(a_1z_1^{\pm 1},\ldots,a_mz_m^{\pm 1};q\bigr)_{\infty}&=
\prod_{j=1}^m\bigl(a_jz_j;q\bigr)_{\infty}\bigl(a_jz_j^{-1};q\bigr)_{\infty}.
\end{split}
\end{equation}
A trigonometric beta integral which closely resembles
\eqref{betaEuler} is the contour integral
\begin{equation}\label{betaEuler1}
\frac{1}{2\pi i}\int_{\mathbb{T}}
\frac{\bigl(-q^{\frac{1}{2}}z^{-1},-q^{\frac{1}{2}}z;q\bigr)_{\infty}}
{\bigl(-q^{\frac{1}{2}}az^{-1},-q^{\frac{1}{2}}bz;q\bigr)_{\infty}}
\frac{dz}{z}=
\frac{\bigl(qa,qb;q\bigr)_{\infty}}{\bigl(q,qab;q\bigr)_{\infty}}
\end{equation}
for $|a|,|b|<|q^{-\frac{1}{2}}|$, where 
\[
\mathbb{T}=\{z\in\mathbb{C} \, | \,  |z|=1 \}
\]
is the positively oriented unit circle in the complex plane. We call 
\eqref{betaEuler1} the Ramanujan integral, since it is closely related
to one of Ramanujan's trigonometric beta integrals from his lost notebook
(see \cite{A2} and references therein).
A second trigonometric analogue 
of \eqref{betaEuler} is the Askey-Wilson integral \cite{AW},
\begin{equation}\label{betaEuler2}
\frac{1}{2\pi i}\int_{\mathbb{T}}
\frac{\bigl(z^{\pm 2};q\bigr)_{\infty}}
{\prod_{j=1}^4 \bigl(t_jz^{\pm 1};q\bigr)_{\infty}}
\frac{dz}{z}=
\frac{2\bigl(t_1t_2t_3t_4;q\bigr)_{\infty}}
{\bigl(q;q\bigr)_{\infty}\prod_{1\leq k<m\leq 4}
\bigl(t_kt_m;q\bigr)_{\infty}}
\end{equation}
for parameters $t_j\in\mathbb{C}$ with $|t_j|<1$ ($j=1,\ldots,4$).
The Askey-Wilson integral is 
the normalization constant for the orthogonality measure of the 
celebrated Askey-Wilson polynomials \cite{AW}.
The Nassrallah-Rahman \cite{NR} integral is the trigonometric beta integral
\begin{equation}\label{betaEuler2n}
\frac{1}{2\pi i}\int_{\mathbb{T}}
\frac{\bigl(z^{\pm 2}, Az^{\pm 1};q\bigr)_{\infty}}
{\prod_{j=0}^4\bigl(t_jz^{\pm 1};q\bigr)_{\infty}}\frac{dz}{z}=
\frac{2\prod_{j=0}^4\bigl(At_j^{-1};q\bigr)_{\infty}}
{\bigl(q;q\bigr)_{\infty}
\prod_{0\leq k<m\leq 4}\bigl(t_kt_m;q\bigr)_{\infty}}
\end{equation}
with $A=t_0t_1t_2t_3t_4$ and with 
parameters $t_j\in\mathbb{C}$ satisfying $|t_j|<1$
($j=0,\ldots,4$). The Nassrallah-Rahman integral
is the normalization constant for the biorthogonality measure of
a five parameter
family of ${}_{10}\phi_9$ rational functions, see \cite{R}.
Note that the Askey-Wilson integral 
\eqref{betaEuler2} is the special case $t_0=0$ of the Nassrallah-Rahman
integral \eqref{betaEuler2n}. The Nassrallah-Rahman integral is the most
general trigonometric analogue of Euler's beta integral \eqref{betaEuler}
known up to date.

One way to introduce 
Ruijsenaars' \cite{R1} hyperbolic gamma function, which will take over the
role of the $q$-Pochhammer symbol in the hyperbolic analogues of 
\eqref{betaEuler1}, \eqref{betaEuler2} and \eqref{betaEuler2n}, 
is by its explicit infinite product realization, see \cite{Sh}, \cite{R2}.
Explicitly, we define for $\tau\in\mathbb{C}\setminus \{0\}$
two deformation parameters $q=q_\tau$ and 
$\widetilde{q}=\widetilde{q}_\tau$ by
\begin{equation}\label{deformationparameters}
q=\exp(2\pi i\tau),\qquad \widetilde{q}=\exp(-2\pi i/\tau).
\end{equation}
The two deformation parameters are thus related by the 
transformation $\tau\mapsto -\frac{1}{\tau}$, which 
is one of the modular 
transformations preserving the upper half plane
\[\mathbb{H}=\{ z\in \mathbb{C} \, | \, \hbox{Im}(z)>0 \}.
\]
We now define for $\tau\in\mathbb{H}$ (so $|q|, |\widetilde{q}|<1$) 
a $(q,\widetilde{q})$-bibasic analogue of the Pochhammer symbol
called the {\it $\tau$-shifted factorial} by
\begin{equation}\label{taufactorial}
\lbrack z;\tau\rbrack_\infty=\frac{\bigl(\exp(-2\pi iz);\widetilde{q}\,
\bigr)_{\infty}}{\bigl(q^z;q\bigr)_{\infty}}.
\end{equation}
Shintani's \cite{Sh}
crucial result implies that $\lbrack z;\tau\rbrack_\infty$ extends 
continuously to $\tau\in\mathbb{R}_{<0}$ for generic
$z\in\mathbb{C}$ (in which 
case $|q|=1=|\widetilde{q}|$). The resulting function, which we still
denote by $\lbrack z;\tau\rbrack_\infty$, 
depends meromorphically on $z\in\mathbb{C}$ and can be expressed explicitly
in terms of Barnes' double gamma function, or equivalently in terms
of Ruijsenaars' \cite{R1} hyperbolic gamma function
or Kurokawa's double sine function (see \cite[Appendix A]{R2}).
In the appendix, \S 6, the interrelation with Ruijsenaars' hyperbolic 
gamma function is made explicit and relevant properties of the hyperbolic
gamma function are translated to the $\tau$-shifted factorial
$\lbrack z;\tau\rbrack_{\infty}$.

The hyperbolic analogues of \eqref{betaEuler1},
\eqref{betaEuler2} and \eqref{betaEuler2n} 
can now be formulated explicitly as follows. 
We define the shorthand notations
\begin{equation*}
\begin{split}
\lbrack a_1,\ldots,a_m;\tau\rbrack_{\infty}&=\prod_{j=1}^m\lbrack a_j;
\tau\rbrack_{\infty},\\
\lbrack a_1\pm z_1,\ldots,a_m\pm z_m;\tau\rbrack_{\infty}&=
\prod_{j=1}^m\lbrack a_j+z_j;\tau\rbrack_{\infty}
\lbrack a_j-z_j;\tau\rbrack_{\infty}
\end{split}
\end{equation*}
for products of $\tau$-shifted factorials.
We use the convention that
$q^u=\exp(2\pi i\tau u)$ and $\widetilde{q}\,{}^u=\exp(-2\pi iu/\tau)$
for $u\in\mathbb{C}$.
The hyperbolic analogue of the Ramanujan integral \eqref{betaEuler1} becomes
\begin{equation}\label{betaEuler3}
\int_{-i\infty}^{i\infty}
\frac{\lbrack \frac{1}{2}+\frac{1}{2\tau}+z,
\frac{1}{2}+\frac{1}{2\tau}-z;\tau\rbrack_{\infty}}
{\lbrack \frac{1}{2}+\frac{1}{2\tau}+\alpha+z,\frac{1}{2}+\frac{1}{2\tau}+
\beta-z;\tau\rbrack_{\infty}}\,dz=
-\frac{q^{-\frac{1}{24}}\widetilde{q}\,{}^{\frac{1}{24}}}{\sqrt{-i\tau}}
\frac{\lbrack \frac{1}{\tau}+
\alpha,\frac{1}{\tau}+\beta;\tau\rbrack_{\infty}}
{\lbrack \frac{1}{\tau}+\alpha+\beta;\tau\rbrack_{\infty}}
\end{equation}
for $\hbox{Re}(\tau\alpha),\hbox{Re}(\tau\beta)<0$ and
$\hbox{Re}\bigl(\alpha-\frac{1}{2}+\frac{1}{2\tau}\bigr),
\hbox{Re}\bigl(\beta-\frac{1}{2}+\frac{1}{2\tau}\bigr)<0$, where
$\sqrt{\cdot}$ is the branch of the square root on 
$\mathbb{C}\setminus \mathbb{R}_{<0}$ which takes positive values on 
$\mathbb{R}_{\geq 0}$ (this branch of the square root will be fixed throughout
the paper). Here the deformation parameter 
$\tau$ may be either from the interior of the second quadrant of the 
complex plane, so $\hbox{Re}(\tau)<0$ and $\hbox{Im}(\tau)>0$ (in which 
case $|q|,|\widetilde{q}|<1$), or from $\mathbb{R}_{<0}$ 
(in which case $|q|=1=|\widetilde{q}|$). The rather mysterious looking
$q$ and $\widetilde{q}$ powers in the right hand side of \eqref{betaEuler3}
arise from an application of the modularity of the Dedekind eta function,
\[\frac{\bigl(q;q\bigr)_{\infty}}
{\bigl(\widetilde{q};\widetilde{q}\,\bigr)_{\infty}}=
\frac{q^{-\frac{1}{24}}\widetilde{q}\,{}^{\frac{1}{24}}}{\sqrt{-i\tau}},
\qquad \tau\in \mathbb{H}.
\]  
The proof of \eqref{betaEuler3}
will given in \S 3. 
With the same conditions on the deformation parameter $\tau$, 
the hyperbolic analogue of the Askey-Wilson integral \eqref{betaEuler2}
becomes
\begin{equation}\label{betaEuler4}
\int_{-i\infty}^{i\infty}
\frac{\lbrack 1\pm 2z;\tau\rbrack_{\infty}}
{\prod_{j=1}^4\lbrack \tau_j\pm z;\tau\rbrack_{\infty}}\,dz
=-\frac{2q^{-\frac{1}{24}}\widetilde{q}\,{}^{\frac{1}{24}}}
{\sqrt{-i\tau}}\,
\frac{\lbrack \tau_1+\tau_2+\tau_3+\tau_3-3;\tau\rbrack_{\infty}}
{\prod_{1\leq k<m\leq 4}\lbrack \tau_k+\tau_m-1;\tau\rbrack_{\infty}}
\end{equation}
for parameters $\tau_j\in \mathbb{C}$ satisfying
$\hbox{Re}(\tau_j)<1$ ($j=1,\ldots,4$) and
$\hbox{Re}\bigl((\tau_1+\tau_2+\tau_3+\tau_4-3)\tau\bigr)<1$. 
The proof will be given in \S 4. 
The hyperbolic analogue of the Nassrallah-Rahman integral
\eqref{betaEuler2n} with deformation parameter $\tau$ 
satisfying the same conditions as for the hyperbolic
Ramanujan integral \eqref{betaEuler3}, reads
\begin{equation}\label{betaEuler5}
\int_{-i\infty}^{i\infty}
\frac{\lbrack 1\pm 2z,a-4\pm z;
\tau\rbrack_{\infty}}
{\prod_{j=0}^4\lbrack \tau_j\pm z;\tau\rbrack_{\infty}}\,dz
=-\frac{2q^{-\frac{1}{24}}\widetilde{q}\,{}^{\frac{1}{24}}}
{\sqrt{-i\tau}}\,
\frac{\prod_{j=0}^4\lbrack a-\tau_j-3;\tau\rbrack_{\infty}}
{\prod_{0\leq k<m\leq 4}\lbrack \tau_k+\tau_m-1;\tau\rbrack_{\infty}}
\end{equation}
where $a=\tau_0+\tau_1+\tau_2+\tau_3+\tau_4$, with the
parameters $\tau_j\in \mathbb{C}$ satisfy
$\hbox{Re}(\tau_j)<1$ ($j=0,\ldots,4$) and
$\hbox{Re}\bigl(a-\tau^{-1} \bigr)>4$. The proof will be given in \S 5.
The hyperbolic Askey-Wilson integral \eqref{betaEuler4} is formally
the limit case $\hbox{Im}(\tau_0)\rightarrow -\infty$
of the hyperbolic Nassrallah-Rahman integral \eqref{betaEuler5},
see Remark \ref{limitNRAW}. 

The proofs are based on 
the observation that a hyperbolic beta integral is essentially
the fusion of a trigonometric sum and a trigonometric integral.
The technique of fusing trigonometric sums and integrals 
will be discussed in \S 2.
The hyperbolic Ramanujan integral \eqref{betaEuler3} is then derived from
fusing Ramanujan's ${}_1\Psi_1$ summation formula with the trigonometric
Ramanujan 
integral \eqref{betaEuler1}, while the hyperbolic Askey-Wilson integral
\eqref{betaEuler4} is derived from fusing Bailey's 
${}_6\Psi_6$ summation formula with the trigonometric Askey-Wilson integral
\eqref{betaEuler2}. The hyperbolic analogue of the Nassrallah-Rahman integral
is more delicate: it is derived from fusing the trigonometric Nassrallah-Rahman
integral with a formula expressing the bilateral very-well-poised
${}_8\Psi_8$ as a sum of two very-well-poised ${}_8\phi_7$ series.
During the fusion procedure the sum of the ${}_8\phi_7$
series is corrected in such a way that it becomes summable by
Bailey's summation formula, leading eventually to \eqref{betaEuler5}.

The fusion technique also reveals a close connection between
trigonometric beta integrals and Macdonald-Mehta type integrals. Roughly
speaking, Macdonald-Mehta integrals can be obtained by fusing 
trigonometric beta integrals
with the inversion formula for the Jacobi theta function
(regarded as a bilateral sum identity). We will encounter several 
one-variable Macdonald-Mehta type integrals in this way.
In particular, we derive
in \S 4 Cherednik's \cite{C} one-variable Macdonald-Mehta integral
as a consequence of the Askey-Wilson integral \eqref{betaEuler2}.

For elliptic beta integrals, the role of Euler's gamma function is
taken over by Ruijsenaars' \cite{R1} elliptic gamma function
\begin{equation}\label{ellgamma}
\Gamma(z;p_1,p_2)=\prod_{k,m=0}^{\infty} \frac{1-z^{-1}p_1^{k+1}p_2^{m+1}}
{1-zp_1^kp_2^m},
\end{equation}
where $p_1,p_2\in\mathbb{C}$ are two arbitrary complex numbers with 
$|p_1|, |p_2|<1$.
Spiridonov \cite{Sp} proved the following elliptic analogue of the
Nassrallah-Rahman integral,
\begin{equation}\label{betaEulerelliptic}
\frac{1}{2\pi i}\int_{\mathbb{T}}
\frac{\prod_{j=0}^4\Gamma(t_jz^{\pm 1};p_1,p_2)}
{\Gamma (z^{\pm 2}, Az^{\pm 1};p_1,p_2)}\,\frac{dz}{z}=
\frac{2}{\bigl(p_1;p_1\bigr)_{\infty}\bigl(p_2;p_2\bigr)_{\infty}}
\frac{\prod_{0\leq k<m\leq 4}\Gamma(t_kt_m;p_1,p_2)}
{\prod_{j=0}^4\Gamma(At_j^{-1};p_1,p_2)}
\end{equation}
where $A=t_0t_1t_2t_3t_4$, with the
five parameters $t_j\in\mathbb{C}$ satisfying
$|t_j|<1$ and $|p_1p_2|<|A|$. Here we used the
notations 
\begin{equation*}
\begin{split}
\Gamma(a_1,\ldots,a_m;p_1,p_2)&=\prod_{j=1}^m\Gamma(a_j;p_1,p_2),\\
\Gamma(a_1z_1^{\pm 1},\ldots,a_mz_m^{\pm 1};p_1,p_2)&=
\prod_{j=1}^m\Gamma(a_jz_j,a_jz_j^{-1};p_1,p_2)
\end{split}
\end{equation*}
for products of elliptic gamma functions. The trigonometric
Nassrallah-Rahman integral \eqref{betaEuler2n} is the special case
of the elliptic Nassrallah-Rahman integral \eqref{betaEulerelliptic}
when one of the deformation parameters $p_j$ is zero.
We show in \S 5.4 that the hyperbolic
Nassrallah-Rahman integral \eqref{betaEuler5} can also be formally
obtained as limit of the elliptic Nassrallah-Rahman 
integral \eqref{betaEulerelliptic}. It is based on a remarkable
limit from the elliptic gamma function to the hyperbolic gamma function
due to Ruijsenaars \cite{R1}. The rigidity of the elliptic theory 
seems to prevent the existence of elliptic degenerations, and consequently 
there do not seem to exist elliptic analogues of 
the Ramanujan integral and of the Askey-Wilson integral.

\begin{rema}
After the appearance of this paper as a preprint,
Erik Koelink kindly pointed out to me that Ruijsenaars 
recently obtained another proof
of the hyperbolic Askey-Wilson integral 
\eqref{betaEuler4} as a spin-off of his studies of  
Hilbert space transforms associated to relativistic hypergeometric functions
in his recent paper \cite{R3}. 
Subsequently van Diejen and Spiridonov 
derived multidimensional analogues of the 
hyperbolic Askey-Wilson and Nassrallah-Rahman
integrals \eqref{betaEuler4} and \eqref{betaEuler5}
by completely different methods in \cite{vDS}. 
Their multidimensional hyperbolic 
Askey-Wilson integral \cite[Thm. 5]{vDS} can be rederived 
with the methods of this paper by fusing Gustafson's \cite{G0} 
multidimensional Askey-Wilson integral with van Diejen's \cite{vD}
multidimensional generalization of Bailey's very-well-poised ${}_6\Psi_6$
summation formula.
\end{rema}

\vspace{.2cm}

{\bf Notations:}
For $|q|<1$ we denote 
\[\bigl(z;q\bigr)_\alpha=\frac{\bigl(z;q\bigr)_{\infty}}
{\bigl(zq^\alpha;q\bigr)_{\infty}}.
\]
This reduces to a finite product when $\alpha\in\mathbb{Z}$.
Products of $\bigl(z;q\bigr)_\alpha$ will be denoted in the same
way as for $\bigl(z;q\bigr)_{\infty}$, see \eqref{productPoch}.
The ${}_{r+1}\phi_{r}$ basic hypergeometric series with base $q$
is defined by
\[{}_{r+1}\phi_{r}\left(\begin{matrix} a_1,a_2,\ldots,a_{r+1}\\
b_1,\ldots,b_r\end{matrix}\,;q,z\right)=
\sum_{m=0}^{\infty}
\frac{\bigl(a_1,a_2,\ldots,a_{r+1};q\bigr)_m}
{\bigl(q,b_1,\ldots,b_r;q\bigr)_m}\,z^m,\qquad |z|<1.
\]
The special case of a very-well-poised ${}_{r+1}\phi_r$ series
is defined as
\begin{equation*}
\begin{split}
{}_{r+1}W_{r}\left(a_1;a_4,\ldots,a_{r+1};q,z\right)&=
{}_{r+1}\phi_{r}
\left(\begin{matrix} a_1,qa_1^{\frac{1}{2}},-qa_1^{\frac{1}{2}},
a_4,\ldots,a_{r+1}\\
a_1^{\frac{1}{2}}, -a_1^{\frac{1}{2}},qa_1/a_4,\ldots,qa_1/a_{r+1}
\end{matrix}\,;q,z\right)\\
&=
\sum_{m=0}^{\infty}\frac{1-a_1q^{2m}}{1-a_1}
\frac{\bigl(a_1,a_4,\ldots,a_{r+1};q\bigr)_m}
{\bigl(q,qa_1/a_4,\ldots,qa_1/a_{r+1};q\bigr)_m}\,z^m.
\end{split}
\end{equation*}
The bilateral basic hypergeometric series ${}_{r}\psi_{r}$ with base
$q$ is defined by
\[{}_{r}\psi_{r}\left(\begin{matrix} a_1,\ldots,a_{r}\\
b_1,\ldots,b_r\end{matrix}\,;q,z\right)=
\sum_{m=-\infty}^{\infty}
\frac{\bigl(a_1,\ldots,a_{r};q\bigr)_m}
{\bigl(b_1,\ldots,b_r;q\bigr)_m}\,z^m,\qquad \left|\frac{b_1\cdots b_r}
{a_1\cdots a_r}\right|<|z|<1.
\]
\vspace{.2cm}\\
{\bf Acknowledgments:}
The author is supported by the Royal Netherlands Academy
of Arts and Sciences (KNAW).



\section{Folding and fusion of integrals.}

In this section we explain the principle of folding and fusion of
integrals, which is fundamental for the study of hyperbolic beta
integrals. 

We take $\tau\in \mathbb{H}$, so that $|q|,|\widetilde{q}\,|<1$ for the
associated deformation parameters $q$ and $\widetilde{q}$ 
(see \eqref{deformationparameters}). 
The (rescaled) Gaussian 
\begin{equation}\label{qGaussian}
G_\tau(z)=q^{\frac{z^2}{2}}
\end{equation}
is analytic, zero-free and 
satisfies the quasi-periodicity conditions
\[G_\tau(z+1)=q^{\frac{1}{2}}q^zG_\tau(z),\qquad
G_\tau\bigl(z-\tau^{-1}\bigr)
=\widetilde{q}\,{}^{-\frac{1}{2}}\exp(-2\pi iz)G_\tau(z).
\]
Suppose that $\phi(z)$ is a 
$\frac{1}{\tau}$-periodic, meromorphic function 
and consider the associated function
\[\widetilde{\phi}(z)=\phi(z)G_\tau(z)=\phi(z)q^{\frac{z^2}{2}}.
\]
By the $\frac{1}{\tau}$-periodicity of $\phi$, we still have
the quasi-periodicity 
\begin{equation}\label{quasi0}
\widetilde{\phi}\bigl(z-\tau^{-1}\bigr)=
\widetilde{q}\,{}^{-\frac{1}{2}}\exp(-2\pi iz)
\widetilde{\phi}(z).
\end{equation}
We now force $\widetilde{\phi}(z)$ to become one-periodic by 
considering the series
\begin{equation}\label{seriesphi}
\phi^+(z)=\sum_{m=-\infty}^{\infty}\widetilde{\phi}(z+m)
\end{equation}
(this resembles the construction of automorphic forms as
Poincar{\'e} series).
We assume that the function $\phi$ is such that the series \eqref{seriesphi}
converges to a meromorphic function. In this situation, the resulting 
bilateral sum $\phi^+(z)$ is thus an {\it one-periodic} meromorphic  
function satisfying
\begin{equation}\label{qper1}
\phi^+\bigl(z-\tau^{-1}\bigr)=
\widetilde{q}\,{}^{-\frac{1}{2}}\exp(-2\pi iz)\phi^+(z).
\end{equation}
On the other hand, the Jacobi theta function
\begin{equation}\label{Jacobi}
\vartheta_{-\frac{1}{\tau}}(z)=
\sum_{m=-\infty}^{\infty}G_{-\frac{1}{\tau}}(m)\exp(2\pi imz)
=\sum_{m=-\infty}^{\infty}\widetilde{q}\,{}^{\frac{m^2}{2}}\exp(2\pi imz)
\end{equation}
is a one-periodic entire function which satisfies the same
quasi-periodicity \eqref{qper1}, hence 
\begin{equation}\label{sumfus}
\phi^+(z)=\sum_{m=-\infty}^{\infty}\widetilde{\phi}(z+m)
=\Phi(z)\vartheta_{-\frac{1}{\tau}}(z)
\end{equation}
with $\Phi$ an elliptic function with respect to the periods $1$ and 
$-\frac{1}{\tau}$. 

If
the bilateral series $\phi^+(z)$ is {\it entire}, then 
the elliptic function $\Phi(z)$ is automatically a constant.
This well known fact follows from the observation 
that the Fourier coefficients 
in the Fourier expansion $\phi^+(z)=\sum_{n}a_n\exp(2\pi inz)$ satisfy
the first order recurrence relation
$a_{n+1}=\widetilde{q}\,{}^{\frac{1}{2}+n}a_n$ ($n\in\mathbb{Z}$)
due to the quasi-periodicity \eqref{qper1} of $\phi^+$.

Choose now in addition a one-periodic meromorphic function $\psi(z)$.
In the examples treated in this paper, the integral
\begin{equation}\label{compactfusion}
\int_0^1\psi(x)\Phi(x)\vartheta_{-\frac{1}{\tau}}(x)dx
\end{equation}
over the period cycle $[0,1]$ will be closely related to 
the trigonometric beta integral \eqref{betaEuler1}, \eqref{betaEuler2}
or \eqref{betaEuler2n}. The integral
\begin{equation}\label{fusion}
\int_{-\infty}^{\infty}\psi(x)\phi(x)q^{\frac{x^2}{2}}dx=
\int_{-\infty}^{\infty}\psi(x)\widetilde{\phi}(x)dx
\end{equation}
may then be viewed as the {\it fusion} of the integral \eqref{compactfusion}
and the sum \eqref{sumfus}, since folding the integral
and using the one-periodicity of $\psi$ implies
\begin{equation}\label{foldsteps}
\begin{split}
\int_{-\infty}^{\infty}\psi(x)\phi(x)q^{\frac{x^2}{2}}dx&=
\sum_{m=-\infty}^{\infty}\int_0^1\psi(x)\widetilde{\phi}(x+m)\,dx\\
&=\int_0^1\psi(x)\phi^+(x)dx=
\int_0^1\psi(x)\Phi(x)\vartheta_{-\frac{1}{\tau}}(x)dx
\end{split}
\end{equation}
by the bilateral sum \eqref{sumfus}. In all applications below, 
\eqref{foldsteps} can be justified by a straightforward application
of Fubini's theorem. Note that the fused integral \eqref{fusion}
admits an explicit evaluation as soon as 
the folded integral \eqref{compactfusion} admits 
some explicit evaluation.

\begin{rema}
A summation formula $F=\sum_{\mathbf{m}\in \mathbb{Z}^n}f(x+\mathbf{m})$
for some function $f:\mathbb{R}^n\rightarrow \mathbb{C}$, with $F$
independent of $x$, implies the identity
\[\iint_{\mathbb{R}^n}f(x)\,dx=F
\]
by the folding technique. 
This observation was extensively used by Gustafson \cite{G} 
to obtain evaluations of multivariate trigonometric integrals from
known multivariate bilateral summation formulas. 
\end{rema}

The simplest example of fusion and folding of integrals arises when
$\phi(z)\equiv 1\equiv \psi(z)$. Then
\[\phi^+(z)=\sum_{m=-\infty}^{\infty}q^{\frac{(z+m)^2}{2}}
=\vartheta_{\tau}(\tau z)q^{\frac{z^2}{2}}
\]
is entire, hence
\begin{equation}\label{sumeasy}
\phi^+(z)=C\vartheta_{-\frac{1}{\tau}}(z)
\end{equation}
is the corresponding sum \eqref{sumfus} for some constant $C$. 
The integral \eqref{compactfusion} is 
\[C\int_0^1\vartheta_{-\frac{1}{\tau}}(x)dx
\]
which equals the constant $C$ in view of 
the explicit Fourier expansion \eqref{Jacobi} 
of the Jacobi theta function. The fused integral \eqref{fusion}
thus becomes
\[\int_{-\infty}^{\infty} q^{\frac{x^2}{2}}dx,
\]
and the above analysis shows that it equals the constant $C$.
The fused integral is up to a change of integration
variable the Gauss integral, whose well known evaluation is given by
\[\int_{-\infty}^{\infty} q^{\frac{x^2}{2}}dx=\frac{1}{\sqrt{-i\tau}}.
\]
We thus conclude that the $C=\frac{1}{\sqrt{-i\tau}}$ in \eqref{sumeasy}. 
Returning to the underlying bilateral sum, we obtain
\begin{equation}\label{inversion}
\sum_{m=-\infty}^{\infty}q^{\frac{(z+m)^2}{2}}=\frac{1}{\sqrt{-i\tau}}
\vartheta_{-\frac{1}{\tau}}(z),
\end{equation}
which is the well known inversion formula for the Jacobi theta function.
The Jacobi inversion formula can be rewritten as
\[\vartheta_{-\frac{1}{\tau}}(z)=
\sqrt{-i\tau}\,q^{\frac{z^2}{2}}\vartheta_\tau(\tau z),
\]
or equivalently, with $z$ replaced by $z/\tau$,
\[
\vartheta_{-\frac{1}{\tau}}(z/\tau)=
\sqrt{-i\tau}\,\widetilde{q}\,{}^{-\frac{z^2}{2}}\vartheta_\tau(z).
\]
Note that in this example, we used the explicit evaluation of the fused
integral to derive the underlying bilateral summation formula 
\eqref{inversion} explicitly.
Conversely, the evaluation of the Gauss integral follows from the 
Jacobi inversion formula and the elementary integral
\begin{equation}\label{thetaintegral}
\int_0^1\vartheta_{-\frac{1}{\tau}}(x)\,dx=1.
\end{equation}

In subsequent sections we show how hyperbolic versions of beta integrals
can be derived by the fusion procedure. In these cases, 
$\psi(x)=\psi_{-\tau^{-1}}(x)$ is a slight modification of the
(one-periodic variant) of the integrand of the trigonometric version 
of the beta integral with respect to base $\widetilde{q}$, while 
$\phi(x)$ is roughly $\psi_\tau(\tau x)^{-1}$ (thus with respect
to the base $q$); this is precisely what one 
would expect from the explicit formula \eqref{taufactorial} 
expressing the hyperbolic gamma function as a quotient of trigonometric
gamma functions.


\section{The hyperbolic Ramanujan integral}

In this section we derive a hyperbolic analogue of the
Ramanujan integral \eqref{betaEuler1}. 
The derivations in this section are based on 
the $q$-binomial theorem
\begin{equation}\label{qbin}
{}_1\phi_0(a;-;q,z)=
\frac{\bigl(az;q\bigr)_{\infty}}{\bigl(z;q\bigr)_{\infty}},\qquad
|q|,|z|<1,
\end{equation}
the Jacobi inversion formula \eqref{inversion}
and the Jacobi triple product identity
\begin{equation}\label{triple}
\vartheta_\tau(z)=
\bigl(q,-q^{\frac{1}{2}}\exp(2\pi iz),-q^{\frac{1}{2}}\exp(-2\pi iz);
q\bigr)_{\infty},\qquad |q|<1.
\end{equation}
An elementary proof of the $q$-binomial theorem can e.g. be found in
\cite{GR}. Of the many known proofs of the Jacobi triple product
identity, we note that there are several which essentially only uses
the $q$-binomium theorem, see e.g. \cite{An}.
We have chosen to give full details of all other identities we encounter
since it clarifies the techniques leading to hyperbolic beta integrals. 


\subsection{Ramanujan's ${}_1\Psi_1$ summation formula}

In this subsection we consider a nontrivial example of the
bilateral sum construction \eqref{sumfus}. 
It leads to an elementary proof of Ramanujan's
${}_1\Psi_1$ summation formula. 
This method of proving bilateral sum identities has
been applied in several closely related setups, see e.g. \cite{Ao}. 

Let $\tau\in \mathbb{H}$ and $0<|a|,|b|<1$. 
We consider the $\frac{1}{\tau}$-periodic entire function
\begin{equation}\label{phi}
\phi(z;a,b)=\phi(z)=\bigl(-aq^{\frac{1}{2}+z},
-bq^{\frac{1}{2}-z};q\bigr)_{\infty}
\end{equation}
and the associated entire function 
$\widetilde{\phi}(z;a,b)=\widetilde{\phi}(z)$ defined by
\[\widetilde{\phi}(z;a,b)=\phi(z;a,b)\,q^{\frac{z^2}{2}}.
\]
The function $\widetilde{\phi}$ satisfies 
\eqref{quasi0}, as well as
$\widetilde{\phi}(z+1)=t(z)\widetilde{\phi}(z)$ 
with 
\[t(z)=\frac{1+b^{-1}q^{\frac{1}{2}+z}}
{1+aq^{\frac{1}{2}+z}}\,b.
\]
As explained in \S 2, we now force $\widetilde{\phi}$ to
become one periodic by considering the bilateral series
\[\phi^+(z)=\sum_{m=-\infty}^{\infty}\widetilde{\phi}(z+m),
\]
which is easily seen to converge absolutely and uniformly on compacta of
$\mathbb{C}$ to an entire function $\phi^+(z)$ due to the conditions on 
$a$ and $b$. Thus we conclude
\[\phi^+(z)=C\,\vartheta_{-\frac{1}{\tau}}(z)
\]
for some constant $C\in\mathbb{C}$. To find an explicit expression for 
the constant $C$,
we first express $\phi^+(z)$ in terms of the ${}_1\Psi_1$ bilateral
series. Using $\widetilde{\phi}(z+1)=
t(z)\widetilde{\phi}(z)$, we can rewrite $\phi^+(z)$ as
\begin{equation}\label{step}
\phi^+(z)=\sum_{m=-\infty}^{\infty}\widetilde{\phi}(z+m)=
t^+(z)\widetilde{\phi}(z)
\end{equation}
with
\begin{equation}\label{s1}
\begin{split}
t^+(z)&=\sum_{m=0}^{\infty}
\prod_{k=0}^{m-1}t(z+k)+
\sum_{m=-\infty}^{-1}\prod_{k=m}^{-1}\frac{1}{t(z+k)}\\
&={}_1\Psi_1\bigl(-b^{-1}q^{\frac{1}{2}+z};
-aq^{\frac{1}{2}+z};q,b\bigr)
\end{split}
\end{equation}
(empty products equal $1$ by convention). Now note that (for generic $a$ and
$b$), $t^+(z)$ can be evaluated at $z\in \mathbb{C}$ satisfying
$q^z=-q^{\frac{1}{2}}a^{-1}$, in which case $t^+(z)$ reduces to a ${}_1\phi_0$
series. By the $q$-binomial theorem we thus conclude
\[t^+(z)|_{q^z=-q^{\frac{1}{2}}a^{-1}}=
\frac{\bigl(q/a;q\bigr)_{\infty}}{\bigl(b;q\bigr)_{\infty}}.
\]
On the other hand,
\[t^+(z)=\frac{\phi^+(z)}{\widetilde{\phi}(z)}=
C\,\frac{\vartheta_{-\frac{1}{\tau}}(z)}{\widetilde{\phi}(z)}=
C\sqrt{-i\tau}\,\frac{\bigl(q,-q^{\frac{1}{2}+z},-q^{\frac{1}{2}-z};
q\bigr)_{\infty}}
{\bigl(-aq^{\frac{1}{2}+z},-bq^{\frac{1}{2}-z};q\bigr)_{\infty}}
\]
where we used the Jacobi inversion formula and the Jacobi triple product
identity for the last equality, hence
\[t^+(z)|_{q^z=-q^{\frac{1}{2}}a^{-1}}=C\sqrt{-i\tau}\,
\frac{\bigl(q/a,a;q\bigr)_{\infty}}{\bigl(ab;q\bigr)_{\infty}}.
\]
We thus conclude that
\begin{equation}\label{C}
C=\frac{1}{\sqrt{-i\tau}}\,\frac{\bigl(ab;q\bigr)_{\infty}}
{\bigl(a,b;q\bigr)_{\infty}}.
\end{equation}
The resulting bilateral summation formula
\begin{equation}\label{Ramanujan}
\phi^+(z;a,b)=
\sum_{m=-\infty}^{\infty}\widetilde{\phi}(z+m;a,b)=
\frac{1}{\sqrt{-i\tau}}
\frac{\bigl(ab;q\bigr)_{\infty}}{\bigl(a,b;q\bigr)_{\infty}}\,
\vartheta_{-\frac{1}{\tau}}(z)
\end{equation}
is Ramanujan's ${}_1\Psi_1$ summation formula,
written in the convenient form \eqref{sumfus}.


\subsection{The Ramanujan integral}
We consider in this subsection the natural trigonometric 
beta integral associated to Ramanujan's ${}_1\Psi_1$ sum \eqref{Ramanujan}.

Again we take $\tau\in\mathbb{H}$, so that $|q|<1$.
Let $\mathbb{T}=\{ z\in \mathbb{C} \, | \, |z|=1\}$ be the positively
oriented unit circle in the complex plane. Choose parameters $c,d\in\mathbb{C}$
satisfying 
$1<|c|<|q^{-\frac{1}{2}}|$ and $0<|d|<|q^{-\frac{1}{2}}|$. We compute the
trigonometric integral
\[\int_{0}^1\frac{\vartheta_\tau(x)}
{\bigl(-q^{\frac{1}{2}}c\exp(-2\pi ix)
,-q^{\frac{1}{2}}d\exp(2\pi ix);q\bigr)_{\infty}}\,dx=
\frac{1}{2\pi i}\int_{\mathbb{T}}\frac{\bigl(q,-q^{\frac{1}{2}}z^{-1},
-q^{\frac{1}{2}}z;q\bigr)_{\infty}}
{\bigl(-q^{\frac{1}{2}}cz^{-1},-q^{\frac{1}{2}}dz;q\bigr)_{\infty}}
\frac{dz}{z}
\]
by shrinking the radius of the integration contour $\mathbb{T}$
to zero while picking up residues (the second expression is obtained 
from the first by the Jacobi triple product identity and
a change of integration variable).
This computation is a simplified version of
a general residual approach to basic contour integrals developed by Slater, 
see \cite[\S 4.9]{GR} for details and references. 

By the conditions
on the parameters, the poles $-cq^{\frac{1}{2}+m}$ ($m\in\mathbb{Z}_{\geq 0}$)
of the integrand lie inside $\mathbb{T}$, while the poles
$-d^{-1}q^{-\frac{1}{2}-m}$ ($m\in\mathbb{Z}_{\geq 0}$) lie outside
$\mathbb{T}$. The condition $|c|>1$ then ensures that
\[\frac{1}{2\pi i}\int_{\mathbb{T}}\frac{\bigl(q,-q^{\frac{1}{2}}z^{-1},
-q^{\frac{1}{2}}z;q\bigr)_{\infty}}
{\bigl(-q^{\frac{1}{2}}cz^{-1},-q^{\frac{1}{2}}dz;q\bigr)_{\infty}}
\frac{dz}{z}=
\sum_{m=0}^{\infty}\,\underset{z=-cq^{\frac{1}{2}+m}}{\hbox{Res}}
\left(\frac{\bigl(q,-q^{\frac{1}{2}}z^{-1},
-q^{\frac{1}{2}}z;q\bigr)_{\infty}}
{\bigl(-q^{\frac{1}{2}}cz^{-1},-q^{\frac{1}{2}}dz;q\bigr)_{\infty}\,z}\right),
\]
cf. \cite[\S 4.9]{GR}. The residues are easily computed, leading to
\begin{equation*}
\begin{split}
\frac{1}{2\pi i}\int_{\mathbb{T}}\frac{\bigl(q,-q^{\frac{1}{2}}z^{-1},
-q^{\frac{1}{2}}z;q\bigr)_{\infty}}
{\bigl(-q^{\frac{1}{2}}cz^{-1},-q^{\frac{1}{2}}dz;q\bigr)_{\infty}}
\frac{dz}{z}&=\frac{\bigl(qc,c^{-1};q\bigr)_{\infty}}
{\bigl(qcd;q\bigr)_{\infty}}\sum_{m=0}^{\infty}
\frac{\bigl(qcd;q\bigr)_m}{\bigl(q;q\bigr)_m}\,c^{-m}\\
&=\frac{\bigl(qc,qd;q\bigr)_{\infty}}{\bigl(qcd;q\bigr)_{\infty}}
\end{split}
\end{equation*}
by the $q$-binomial theorem. The conditions $|c|>1$ and $|d|>0$
may now be removed by analytic continuation.
This proves the Ramanujan integral 
\eqref{betaEuler1},
\begin{equation}\label{qbetasimple}
\int_0^1\frac{\vartheta_\tau(x)}
{\bigl(-q^{\frac{1}{2}}c\exp(-2\pi ix),-q^{\frac{1}{2}}d\exp(2\pi ix);
q\bigr)_{\infty}}\,dx=
\frac{\bigl(qc,qd;q\bigr)_{\infty}}{\bigl(qcd;q\bigr)_{\infty}}
\end{equation}
for  $\tau\in\mathbb{H}$ and $|c|,|d|<|q^{-\frac{1}{2}}|$.
The connection with the integrals from
Ramanujan's lost notebook will become apparent in the following subsection.
\subsection{Fusion}
The Ramanujan ${}_1\Psi_1$ summation formula \eqref{Ramanujan}
and the trigonometric Ramanujan integral \eqref{qbetasimple} can be 
fused as follows.
\begin{prop}\label{fusionRam}
Let $\tau\in\mathbb{H}$, $0<|a|,|b|<1$ and 
$|c|,|d|<|\widetilde{q}\,{}^{-\frac{1}{2}}|$. Then
\begin{equation}\label{fusion1}
\int_{-\infty}^{\infty}\frac{\bigl(-aq^{\frac{1}{2}+x},
-bq^{\frac{1}{2}-x};q\bigr)_{\infty}\,q^{\frac{x^2}{2}}}
{\bigl(-\widetilde{q}^{\frac{1}{2}}c\exp(-2\pi ix),
-\widetilde{q}^{\frac{1}{2}}d\exp(2\pi ix);\widetilde{q}\,\bigr)_{\infty}}
\,dx=\frac{1}{\sqrt{-i\tau}}
\frac{\bigl(ab;q\bigr)_{\infty}}{\bigl(\widetilde{q}cd;\widetilde{q}\,
\bigr)_{\infty}}\frac{\bigl(\widetilde{q}c,
\widetilde{q}d;\widetilde{q}\,\bigr)_{\infty}}
{\bigl(a,b;q\bigr)_{\infty}}.
\end{equation}
\end{prop}
\begin{proof}
The integrand $J(x)$ of \eqref{fusion1} can be written as
\[J(x)=\frac{\widetilde{\phi}(x;a,b)}
{\bigl(-\widetilde{q}^{\frac{1}{2}}c\exp(-2\pi ix),
-\widetilde{q}^{\frac{1}{2}}d\exp(2\pi ix);
\widetilde{q}\,\bigr)_{\infty}}.
\]
Note that the denominator is one-periodic.
Hence folding the integral gives
\begin{equation*}
\begin{split}
\int_{-\infty}^{\infty}J(x)\,dx&=
\int_0^1\frac{\phi^+(x;a,b)}
{\bigl(-\widetilde{q}^{\frac{1}{2}}c\exp(-2\pi ix),
-\widetilde{q}^{\frac{1}{2}}d\exp(2\pi ix);
\widetilde{q}\,\bigr)_{\infty}}\,dx\\
&=\frac{1}{\sqrt{-i\tau}}\frac{\bigl(ab;q\bigr)_{\infty}}
{\bigl(a,b;q\bigr)_{\infty}}
\int_0^1\frac{\vartheta_{-\frac{1}{\tau}}(x)}
{\bigl(-\widetilde{q}^{\frac{1}{2}}c\exp(-2\pi ix),
-\widetilde{q}^{\frac{1}{2}}d\exp(2\pi ix);
\widetilde{q}\,\bigr)_{\infty}}\,dx\\
&=\frac{1}{\sqrt{-i\tau}}\frac{\bigl(ab;q\bigr)_{\infty}}
{\bigl(\widetilde{q}cd;\widetilde{q}\,\bigr)_{\infty}}
\frac{\bigl(\widetilde{q}c,\widetilde{q}d;\widetilde{q}\,\bigr)_{\infty}}
{\bigl(a,b;q\bigr)_{\infty}},
\end{split}
\end{equation*}
where we used Ramanujan's ${}_1\Psi_1$ summation formula \eqref{Ramanujan}
for the second equality
and the trigonometric beta integral \eqref{qbetasimple} for the last identity.
\end{proof}
The proposition is also valid for $a=b=0$, in which case the resulting integral
\begin{equation}\label{Ram1n}
\int_{-\infty}^{\infty}
\frac{q^{\frac{x^2}{2}}}
{\bigl(-\widetilde{q}^{\frac{1}{2}}c\exp(-2\pi ix),
-\widetilde{q}^{\frac{1}{2}}d\exp(2\pi ix);\widetilde{q}\,\bigr)_{\infty}}
\,dx=\frac{1}{\sqrt{-i\tau}}
\frac{\bigl(\widetilde{q}c,\widetilde{q}d;\widetilde{q}\,\bigr)_{\infty}}
{\bigl(\widetilde{q}cd;\widetilde{q}\,\bigr)_{\infty}}
\end{equation}
is the fusion of  
the Jacobi inversion formula \eqref{inversion} and the trigonometric
Ramanujan integral \eqref{qbetasimple}. 
The integral \eqref{Ram1n} is one of Ramanujan's trigonometric
analogues of the beta integral from his lost notebook, see \cite{A2},
\cite[Exerc. 6.15]{GR} and references therein. Ramanujan's second
trigonometric analogue of the beta integral
is the special case $c=d=0$ of the proposition,
\begin{equation}\label{Ram2n}
\int_{-\infty}^{\infty}\bigl(-aq^{\frac{1}{2}+x},
-bq^{\frac{1}{2}-x};q\bigr)_{\infty}\,q^{\frac{x^2}{2}}\,dx=
\frac{1}{\sqrt{-i\tau}}
\frac{\bigl(ab;q\bigr)_{\infty}}{\bigl(a,b;q\bigr)_{\infty}}.
\end{equation}
This integral is the fusion of 
Ramanujan's ${}_1\Psi_1$ summation formula \eqref{Ramanujan}
and the elementary integral \eqref{thetaintegral}.

\subsection{The hyperbolic Ramanujan integral}

The integral
\eqref{fusion1} turns out to be more than just a
particular $(q,\widetilde{q}\,)$-bibasic extension of the Ramanujan 
integral. As we will see in this subsection, it is the crucial 
intermediate step towards the hyperbolic Ramanujan integral.

We specialize the parameters in the fused Ramanujan integral
\eqref{fusion1} to
\begin{equation}\label{parspecialRam}
a=q^\alpha,\qquad b=q^\beta,\qquad c=\exp(-2\pi i\alpha),\qquad
d=\exp(-2\pi i\beta).
\end{equation}
We assume that $\tau\in \mathbb{H}$ and that $\alpha,\beta\in \mathbb{C}$
are such that the corresponding
parameters $a,b,c,d$ satisfy the parameter constraints of 
Proposition \ref{fusionRam}. The fused Ramanujan integral 
\eqref{fusion1} can then be written entirely in terms of the 
$\tau$-shifted factorial \eqref{taufactorial},
\begin{equation}\label{step2R}
\int_{-\infty}^{\infty}
\frac{q^{\frac{x^2}{2}}}{\lbrack \frac{1}{2}+\frac{1}{2\tau}+\alpha+x,
\frac{1}{2}+\frac{1}{2\tau}+\beta-x;\tau\rbrack_\infty}\,dx=
\frac{1}{\sqrt{-i\tau}}\frac{\lbrack \frac{1}{\tau}+\alpha,
\frac{1}{\tau}+\beta;\tau\rbrack_\infty}
{\lbrack \frac{1}{\tau}+\alpha+\beta;\tau\rbrack_{\infty}}.
\end{equation}
By the reflection equation \eqref{reflectioneq} we can rewrite 
\eqref{step2R} in a form which closely resembles the trigonometric
Ramanujan integral \eqref{betaEuler1},
\begin{equation}\label{step3R}
\int_{-\infty}^{\infty}
\frac{\lbrack \frac{1}{2}+\frac{1}{2\tau}+x,\frac{1}{2}+\frac{1}{2\tau}-x;
\tau\rbrack_{\infty}}{\lbrack \frac{1}{2}+\frac{1}{2\tau}+\alpha+x,
\frac{1}{2}+\frac{1}{2\tau}+\beta-x;\tau\rbrack_\infty}\,dx=
\frac{q^{-\frac{1}{24}}\widetilde{q}\,{}^{\frac{1}{24}}}
{\sqrt{-i\tau}}\frac{\lbrack \frac{1}{\tau}+\alpha,
\frac{1}{\tau}+\beta;\tau\rbrack_\infty}
{\lbrack \frac{1}{\tau}+\alpha+\beta;\tau\rbrack_{\infty}}.
\end{equation}
One may view \eqref{step3R} already as an hyperbolic analogue
of the Ramanujan integral. It is though not completely satisfactory, 
because it does not extend to 
$\tau\in \mathbb{R}_{<0}$ (in which case the deformation parameters
$q$ and $\widetilde{q}$ satisfy $|q|=|\widetilde{q}\,|=1$).
To include this important case, we first need to rotate the integration cycle
$\mathbb{R}$ to $-i\mathbb{R}$ about the origin. 

For the rotation of the integration cycle it 
is convenient to restrict attention to $\tau\in \mathbb{H}\cap\mathbb{C}_-$,
where $\mathbb{C}_-$ (respectively $\mathbb{C}_+$) 
is the open left (respectively right) half plane in $\mathbb{C}$, 
and to parameters $\alpha,\beta\in \mathbb{C}$ satisfying
\begin{equation}\label{pardomainRam}
0<\alpha,\beta<\frac{1}{2}.
\end{equation}
With these parameter constraints,
the corresponding parameters $a,b,c,d$ (see \eqref{parspecialRam})
automatically satisfy the conditions of Proposition \ref{fusionRam}.
We write
\[I(z)=\frac{\lbrack \frac{1}{2}+\frac{1}{2\tau}+z,
\frac{1}{2}+\frac{1}{2\tau}-z;
\tau\rbrack_{\infty}}{\lbrack \frac{1}{2}+\frac{1}{2\tau}+\alpha+z,
\frac{1}{2}+\frac{1}{2\tau}+\beta-z;\tau\rbrack_\infty}
\]
for the integrand of \eqref{step3R}. By the zero and pole
locations of the $\tau$-shifted factorial (see \eqref{zero}
and \eqref{pole} respectively)
we conclude that the poles of $I(z)$ are contained in the two sets
\begin{equation}\label{sets}
\begin{split}
-&\alpha+\frac{1}{2}-\frac{1}{2\tau}+\frac{1}{\tau}\mathbb{Z}_{\leq 0}+
\mathbb{Z}_{\geq 0},\\
&\beta-\frac{1}{2}+\frac{1}{2\tau}+\frac{1}{\tau}\mathbb{Z}_{\geq 0}+
\mathbb{Z}_{\leq 0}.
\end{split}
\end{equation}
By the conditions $\tau\in \mathbb{H}\cap \mathbb{C}_-$ and
\eqref{pardomainRam}, the first sequence
$-\alpha+\frac{1}{2}-\frac{1}{2\tau}+\frac{1}{\tau}\mathbb{Z}_{\leq 0}+
\mathbb{Z}_{\geq 0}$ lies in the interior of the 
first quadrant of the complex plane (so in $\mathbb{H}\cap\mathbb{C}_+$),
while the second sequence
$\beta-\frac{1}{2}+\frac{1}{2\tau}+\frac{1}{\tau}\mathbb{Z}_{\geq 0}+
\mathbb{Z}_{\leq 0}$
lies in the interior of the third quadrant of the complex plane
(so in $(-\mathbb{H})\cap\mathbb{C}_-$).
Hence the integrand $I(z)$ of \eqref{step3R}
is analytic in the
second and fourth quadrant of the complex plane. 
We now clockwise rotate the integration contour 
$\mathbb{R}$ to $-i\mathbb{R}$ in \eqref{step3R} about the origin. We thus 
stay within the union of the second and fourth quadrant 
and do not pass poles of $I(z)$. 
To rigorously verify that the integral evaluation does not
alter, we need to take the asymptotic behaviour of 
the integrand $I(z)$ into account.

\begin{lem}\label{boundRam}
Suppose that $\tau\in \mathbb{H}\cap \mathbb{C}_-$ and that
the parameters satisfy \eqref{pardomainRam}. 
Then there exists a constant $C\in \mathbb{R}_{>0}$ such that
\begin{equation*}
\begin{split}
|I(z)|&\leq C|q^{\beta z}|\,\,\,\,\, \hbox{ if }\,\,
\hbox{Re}(z)\geq 0\,\,\hbox{ and }\,\, \hbox{Im}(z)\leq 0,\\
|I(z)|&\leq C|q^{-\alpha z}|\,\, \hbox{ if }\,\,
\hbox{Re}(z)\leq 0\,\,\hbox{ and }\,\, \hbox{Im}(z)\geq 0.
\end{split}
\end{equation*}
\end{lem}
\begin{proof}
Using the reflection equation \eqref{reflectioneq}
we can write $I(z)$ as
\begin{equation}\label{rewrittenI}
\begin{split}
I(z)&=
C_1q^{\beta z}\,\frac{\lbrack \frac{1}{2}+\frac{1}{2\tau}-\beta+z;
\tau\rbrack_{\infty}}
{\lbrack \frac{1}{2}+\frac{1}{2\tau}+\alpha+z;\tau\rbrack_{\infty}}\\
&=C_2q^{-\alpha z}\,\frac{\lbrack \frac{1}{2}+\frac{1}{2\tau}-\alpha-z;\tau
\rbrack_{\infty}}
{\lbrack \frac{1}{2}+\frac{1}{2\tau}+\beta-z;\tau\rbrack_{\infty}},
\end{split}
\end{equation}
for certain irrelevant $z$-independent nonzero constants $C_1$ and $C_2$.
For the asymptotics in the fourth quadrant we use now the first expression
for $I(z)$ in \eqref{rewrittenI}. 
Using the infinite product expression \eqref{taufactorial}
for the $\tau$-shifted factorial, we may thus write
\begin{equation}\label{rewrittenI2}
I(z)=C_1\frac{\bigl(-\widetilde{q}\,{}^{\frac{1}{2}}d^{-1}\exp(-2\pi iz);
\widetilde{q}\,\bigr)_{\infty}}
{\bigl(-\widetilde{q}\,{}^{\frac{1}{2}}c\exp(-2\pi iz);
\widetilde{q}\,\bigr)_{\infty}}
\frac{\bigl(-q^{\frac{1}{2}+\alpha}q^z;q\bigr)_{\infty}}
{\bigl(-q^{\frac{1}{2}-\beta}q^z;q\bigr)_{\infty}}\,q^{\beta z}
\end{equation}
with $c$ and $d$ given by \eqref{parspecialRam}. 
Now note that $z\mapsto \exp(-2\pi iz)$
maps the fourth quadrant $\{z\in \mathbb{C} \,|\, 
\hbox{Re}(z)\geq 0,\,\, \hbox{Im}(z)\leq 0\}$ of the complex
plane into the closed
unit disc 
\[
\mathbb{D}=\{w\in \mathbb{C} \, | \, |w|\leq 1\}.
\] 
Since furthermore $|\widetilde{q}\,{}^{\frac{1}{2}}|<1$ and
$|c|=1$ by the parameter constraints
\eqref{pardomainRam}, we conclude that
\[\left|
\frac{\bigl(-\widetilde{q}\,{}^{\frac{1}{2}}d^{-1}\exp(-2\pi iz);
\widetilde{q}\,\bigr)_{\infty}}
{\bigl(-\widetilde{q}\,{}^{\frac{1}{2}}c\exp(-2\pi iz);
\widetilde{q}\,\bigr)_{\infty}}\right|
\]
is uniformly bounded on the fourth quadrant of the complex plane.
Observe that $z\mapsto q^z$ maps the fourth quadrant
of the complex plane into the closed unit disc $\mathbb{D}$
since $\tau\in \mathbb{H}\cap \mathbb{C}_-$. 
Furthermore, $|q^{\frac{1}{2}-\beta}|<1$
by the parameter constraints \eqref{pardomainRam}, hence 
\[
\left|\frac{\bigl(-q^{\frac{1}{2}+\alpha}q^z;q\bigr)_{\infty}}
{\bigl(-q^{\frac{1}{2}-\beta}q^z;q\bigr)_{\infty}}\right|
\]
is uniformly bounded on the fourth quadrant of the complex plane.
In view of \eqref{rewrittenI2} this proves 
the uniform asymptotics for $I(z)$ with $z$ in the
fourth quadrant of the complex plane. 
The asymptotics of $I(z)$ in the second quadrant of the complex
plane is determined similarly
using the second expression for $I(z)$
in \eqref{rewrittenI}. 
\end{proof}
We have the following direct consequence of Lemma \ref{boundRam}.
\begin{cor}\label{boundRamcor}
Let $\tau\in \mathbb{H}\cap \mathbb{C}_-$ and suppose that the parameters 
$\alpha,\beta\in\mathbb{C}$ satisfy \eqref{pardomainRam}. 
Let $0<\epsilon<1$ be the growth coefficient
\[\epsilon=\exp(-2\pi \min(\alpha,\beta)\,m)\]
with $m$ the strictly positive constant 
\[m=\min_{\theta\in [-\frac{\pi}{2},0]}
\left(\hbox{Im}(\tau e^{i\theta})\right).
\]
Then there exists a constant $C\in \mathbb{R}_{>0}$ 
such that
\[|I(z)|\leq C\epsilon^{|z|}
\]
for $z\in \mathbb{C}$ satisfying  $\hbox{Re}(z)\geq 0$ and 
$\hbox{Im}(z)\leq 0$, as well as for $z\in\mathbb{C}$ satisfying
$\hbox{Re}(z)\leq 0$ and $\hbox{Im}(z)\geq 0$.
\end{cor}

We keep the assumptions $\tau\in \mathbb{H}\cap \mathbb{C}_-$
and \eqref{pardomainRam}. In view of Corollary \ref{boundRamcor},
we may apply Cauchy's theorem to rotate clockwise the integration
contour $\mathbb{R}$ in \eqref{step3R} to $-i\mathbb{R}$ about the origin
without altering its evaluation.
Relaxing the parameter constraints
leads now to the following main result.
\begin{thm}[Hyperbolic Ramanujan integral]\label{hyperbolicR}
Let $\tau\in \mathbb{C}$ with $\hbox{Re}(\tau)<0$ and 
$\hbox{Im}(\tau)\geq 0$, and suppose that the parameters
$\alpha,\beta\in\mathbb{C}$ satisfy
\begin{equation}\label{pardomainRam2}
\hbox{Re}(\tau\alpha), \hbox{Re}(\tau\beta)<0,
\qquad \hbox{Re}\bigl(\alpha-\frac{1}{2}+\frac{1}{2\tau}\bigr),
\hbox{Re}\bigl(\beta-\frac{1}{2}+\frac{1}{2\tau}\bigr)<0.
\end{equation}
Then
\begin{equation}\label{step6R}
\int_{-i\infty}^{i\infty}
\frac{\lbrack \frac{1}{2}+\frac{1}{2\tau}+z,\frac{1}{2}+\frac{1}{2\tau}-z;
\tau\rbrack_{\infty}}{\lbrack \frac{1}{2}+\frac{1}{2\tau}+\alpha+z,
\frac{1}{2}+\frac{1}{2\tau}+\beta-z;\tau\rbrack_\infty}\,dz=
-\frac{q^{-\frac{1}{24}}\widetilde{q}\,{}^{\frac{1}{24}}}
{\sqrt{-i\tau}}\frac{\lbrack \frac{1}{\tau}+\alpha,
\frac{1}{\tau}+\beta;\tau\rbrack_\infty}
{\lbrack \frac{1}{\tau}+\alpha+\beta;\tau\rbrack_{\infty}}.
\end{equation}
\end{thm}
\begin{proof}
Let $\tau\in \mathbb{H}\cap \mathbb{C}_-$.
We already argued that \eqref{step6R} is valid for parameters
$\alpha$ and $\beta$ satisfying the constraints \eqref{pardomainRam}.
By analytic continuation one easily verifies that the integral evaluation
\eqref{step6R} is valid under the milder constraints \eqref{pardomainRam2}.
In fact, in view of \eqref{as1} and \eqref{as2} 
the requirements that the integrand $I(z)$ decays 
exponentially for $z\rightarrow \pm i\infty$ and that the sequence of poles 
$-\alpha+\frac{1}{2}-\frac{1}{2\tau}+\frac{1}{\tau}\mathbb{Z}_{\leq 0}
+\mathbb{Z}_{\geq 0}$ (respectively 
$\beta-\frac{1}{2}+\frac{1}{2\tau}+\frac{1}{\tau}\mathbb{Z}_{\geq 0}+
\mathbb{Z}_{\leq 0}$) of $I(z)$ is contained in $\mathbb{C}_+$
(respectively $\mathbb{C}_-$), lead to the conditions \eqref{pardomainRam2}
on the parameters $\alpha$ and $\beta$. 

Let $\tau^\prime\in \mathbb{R}_{<0}$ and choose 
$\alpha,\beta\in \mathbb{C}$ satisfying \eqref{pardomainRam2}
with $\tau=\tau^\prime$. Let $0<\delta<\frac{\pi}{2}$ and define the 
compact subset
\[K=\{\tau(\phi)\, | \, \phi\in [0,\delta] \}\subset \mathbb{C}_-
\]
with $\tau(\phi)=\tau^\prime\exp(-i\phi)$.
By choosing $\delta$ small enough we may assume
that the parameter constraints \eqref{pardomainRam2} are valid
for all $\tau\in K$. Observe that 
$K\cap \mathbb{H}=\{\tau(\phi) \, | \, 0<\phi\leq \delta\}$
and that $K\cap \mathbb{R}_{<0}=\tau(0)=\tau^\prime$. By 
the reflection equation \eqref{reflectioneq} and by the asymptotic estimates
\eqref{as1} and \eqref{as2} for the $\tau$-shifted factorial, the weight 
function $I(z)$ in \eqref{step6R} has the asymptotic behaviour
\begin{equation*}
I(z)=
\begin{cases}
\mathcal{O}\bigl(\exp(-2\pi i\tau\alpha z)\bigr),\qquad &\hbox{ if } 
\hbox{Im}(z)\rightarrow\infty,\\
\mathcal{O}\bigl(\exp(2\pi i\tau\beta z)\bigr),\qquad &\hbox{ if }
\hbox{Im}(z)\rightarrow -\infty,
\end{cases}
\end{equation*}
uniformly for $\tau\in K$ and for $\hbox{Re}(z)$ in compacta
of $\mathbb{R}$. Consequently,
\[|I(ix)|=\mathcal{O}(\epsilon^{|x|}),\qquad x\mapsto \pm\infty\]
uniformly for $\tau\in K$, with growth exponent
\[\epsilon=\exp\left(2\pi \max_{\tau\in K}\bigl(\hbox{Re}(\tau\alpha),
\hbox{Re}(\tau\beta)\bigr)\right)<1
\]
depending only on $\alpha,\beta$ and $K$. Thus if we take
$\tau=\tau(\phi)=K\cap\mathbb{H}$ in \eqref{step6R},
then Lebesgue's dominated convergence theorem allows us to
interchange the limit $\phi\searrow 0$ with the integral
over $i\mathbb{R}$. This
yields \eqref{step6R} for $\tau=\tau^\prime\in \mathbb{R}_{<0}$.
\end{proof}


\section{The hyperbolic Askey-Wilson integral}

In this section we apply the folding and fusion 
techniques to derive an hyperbolic analogue of the Askey-Wilson
integral. For proofs of the underlying 
summation formula and the underlying trigonometric Askey-Wilson integral
\eqref{betaEuler2} we refer now directly to the literature 
(in practice we often refer to the book 
\cite{GR} of Gasper and Rahman, in which one can find further
detailed references to the literature).

\subsection{Bailey's ${}_6\Psi_6$ summation formula}

We take $\tau\in\mathbb{H}$, so that $|q|<1$.
For the Askey-Wilson integral and its hyperbolic analogue, the role
of the function $\phi$ in the underlying summation formula
(see \S 2) is 
\begin{equation}\label{phi2}
\phi(z)=
\frac{\prod_{j=1}^4\bigl(s_jq^{\pm z};q\bigr)_{\infty}}
{\bigl(q^{1\pm z}, -q^{1\pm z},   
q^{\frac{1}{2}\pm z};q\bigr)_{\infty}}.
\end{equation}
As in \S 2, we define 
\[
\widetilde{\phi}(z)=\phi(z)\,q^{\frac{z^2}{2}}.
\]
We assume for the remainder of this subsection that the parameters $s_j$ 
are nonzero complex parameters and satisfy 
$|q^{-3}s_1s_2s_3s_4|<1$. A direct computation shows that
\[\widetilde{\phi}(z+1)=t(z)\widetilde{\phi}(z)
\]
with 
\[t(z)=\frac{s_1s_2s_3s_4}{q^3}\frac{\bigl(1-q^{1+z}\bigr)
\bigl(1+q^{1+z}\bigr)}{\bigl(1-q^z\bigr)\bigl(1+q^z\bigr)}
\prod_{j=1}^4\frac{\bigl(1-s_j^{-1}q^{1+z}\bigr)}
{\bigl(1-s_jq^z\bigr)}.
\]
Due to the conditions on the parameters the bilateral sum 
\begin{equation*}
\begin{split}
t^+(z)=&\sum_{m=0}^{\infty}\prod_{k=0}^{m-1}t(z+k)+
\sum_{m=-\infty}^{-1}\prod_{k=m}^{-1}\frac{1}{t(z+k)}\\
=&{}_6\Psi_6\left(\begin{matrix} q^{1+z}, -q^{1+z}, s_1^{-1}q^{1+z},
s_2^{-1}q^{1+z}, s_3^{-1}q^{1+z}, s_4^{-1}q^{1+z}\\
q^z, -q^z, s_1q^z, s_2q^z, s_3q^z, s_4q^z\end{matrix};q, 
\frac{s_1s_2s_3s_4}{q^3}\right)
\end{split}
\end{equation*}
converges absolutely and uniformly on compacta away from the 
$\mathbb{Z}_{\leq 0}$-translates of the poles and the
$\mathbb{Z}_{>0}$-translates of the zeros of $t(z)$, and it
extends to a meromorphic function in $z\in\mathbb{C}$.
Hence the bilateral sum
\[\phi^+(z)=\sum_{m=-\infty}^{\infty}\widetilde{\phi}(z+m)
\]
converges to a meromorphic function in $z\in\mathbb{C}$, given explicitly 
by
\[\phi^+(z)=t^+(z)\widetilde{\phi}(z).
\]
We now take a shortcut compared to our analysis of
the Ramanujan integrals by directly applying Bailey's ${}_6\Psi_6$
summation formula \cite[(5.3.1)]{GR}
for the explicit evaluation of $t^+(z)$,
\[
t^+(z)=
\frac{\bigl(q;q\bigr)_{\infty}\prod_{1\leq k<m\leq 4}
\bigl(q^{-1}s_ks_m;q\bigr)_{\infty}}
{\bigl(q^{-3}s_1s_2s_3s_4;q\bigr)_{\infty}}
\frac{\bigl(q^{1\pm 2z};q\bigr)_{\infty}}
{\prod_{j=1}^4\bigl(s_jq^{\pm z};q\bigr)_{\infty}}.
\]
Since
\begin{equation}\label{expanding}
\bigl(u^2;q\bigr)_{\infty}=\bigl(u,-u,q^{\frac{1}{2}}u,
-q^{\frac{1}{2}}u;q\bigr)_{\infty},
\end{equation}
the meromorphic function $\phi^+(z)=t^+(z)\widetilde{\phi}(z)$ 
can now be written as
\begin{equation}\label{Bailey}
\phi^+(z)=\frac{1}{\sqrt{-i\tau}}
\frac{\prod_{1\leq k<m\leq 4}\bigl(q^{-1}s_ks_m;q\bigr)_{\infty}}
{\bigl(q^{-3}s_1s_2s_3s_4;q\bigr)_{\infty}}\,
\vartheta_{-\frac{1}{\tau}}(z)
\end{equation}
after a straightforward computation using the Jacobi triple product identity
and the Jacobi inversion formula. Written in the form \eqref{sumfus},
we thus have $\phi^+(z)=\Phi(z)\vartheta_{-\frac{1}{\tau}}(z)$ with
the elliptic function $\Phi(z)$ being the constant
\begin{equation}\label{PhiAW}
\Phi=\frac{1}{\sqrt{-i\tau}}
\frac{\prod_{1\leq k<m\leq 4}
\bigl(q^{-1}s_ks_m;q\bigr)_{\infty}}
{\bigl(q^{-3}s_1s_2s_3s_4;q\bigr)_{\infty}}.
\end{equation}

\subsection{Fusion}
In this subsection we fuse the trigonometric Askey-Wilson integral
\eqref{betaEuler2} with Bailey's ${}_6\Psi_6$ summation formula 
\eqref{Bailey}. For proofs of the trigonometric Askey-Wilson integral 
\eqref{betaEuler2}, see \cite{AW}, \cite[\S 6]{GR} and references therein.
\begin{prop}\label{fusedAW}
Let $\tau\in \mathbb{H}$ and $s_j,t_k\in \mathbb{C}$ 
\textup{(}$j,k=1,\ldots,4$\textup{)} with $|q^{-3}s_1s_2s_3s_4|<1$ and
$|t_k|<1$ \textup{(}$k=1,\ldots,4$\textup{)}. Then
\begin{equation*}
\begin{split}
\int_{-\infty}^{\infty}
\frac{\bigl(\exp(\pm 2\pi ix),-\exp(\pm 2\pi ix), 
\widetilde{q}\,{}^{\frac{1}{2}}\exp(\pm 2\pi ix);
\widetilde{q}\,\bigr)_{\infty}}
{\bigl(q^{1\pm x},-q^{1\pm x}, q^{\frac{1}{2}\pm x};q\bigr)_{\infty}}
&\prod_{j=1}^4\frac{\bigl(s_jq^{\pm x};q\bigr)_{\infty}}
{\bigl(t_j\exp(\pm 2\pi ix);\widetilde{q}\,\bigr)_{\infty}}\,
q^{\frac{x^2}{2}}dx\\
=\frac{2}{\sqrt{-i\tau}}
\frac{\bigl(t_1t_2t_3t_4;\widetilde{q}\,\bigr)_{\infty}}
{\bigl(q^{-3}s_1s_2s_3s_4;q\bigr)_{\infty}}
&\prod_{1\leq k<m\leq 4}\frac{\bigl(q^{-1}s_ks_m;q\bigr)_{\infty}}
{\bigl(t_kt_m;\widetilde{q}\,\bigr)_{\infty}}.
\end{split}
\end{equation*}
\end{prop}
\begin{proof}
The integrand $J(x)$ can be written as
\[J(x)=\widetilde{\phi}(x)w(x)
\]
with $\widetilde{\phi}(x)$
as in \S 4.1 and with $w(x)$
the one-periodic function
\[w(x)=
\frac{\bigl(\exp(\pm 2\pi ix),-\exp(\pm 2\pi ix), 
\widetilde{q}\,{}^{\frac{1}{2}}
\exp(\pm 2\pi ix);\widetilde{q}\,\bigr)_{\infty}}
{\prod_{j=1}^4\bigl(t_j\exp(\pm 2\pi ix);\widetilde{q}\,\bigr)_{\infty}}.
\]
Note that $J(x)$ is regular on $\mathbb{R}$ since the
real zeros of the factor
$\bigl(q^{1\pm x}, q^{\frac{1}{2}\pm x};q\bigr)_{\infty}$
in the denominator of $J(x)$ are compensated by zeros of the factor 
$\bigl(\exp(\pm 2\pi ix),-\exp(\pm 2\pi ix);\widetilde{q}\,\bigr)_{\infty}$
in the numerator of $J(x)$. 
Thus folding the integral and using Bailey's ${}_6\Psi_6$ sum 
\eqref{Bailey}, we obtain with $\Phi$ the constant \eqref{PhiAW},
\begin{equation*}
\begin{split}
\int_{-\infty}^{\infty}J(x)\,dx&=
\int_{-\infty}^{\infty}\widetilde{\phi}(x)w(x)\,dx\\
&=\int_0^1\phi^+(x)w(x)\,dx
=\Phi \int_0^1\vartheta_{-\frac{1}{\tau}}(x)w(x)\,dx
\end{split}
\end{equation*}
by Fubini's theorem. Using now the Jacobi triple product identity
and \eqref{expanding}, 
the resulting
integral over the period cycle $[0,1]$ is essentially
the trigonometric Askey-Wilson
integral \eqref{betaEuler2} in base $\widetilde{q}$. Its evaluation 
gives the desired result.  
\end{proof}
\begin{rema}\label{polAWremark}
Under the parameter constraints of Proposition \ref{fusedAW} we have
\[\int_{-\infty}^{\infty}f(x)J(x)\,dx=
\Phi\int_{0}^1f(x)\vartheta_{-\frac{1}{\tau}}(x)w(x)\,dx
\]
for a regular one-period function $f:\mathbb{R}\rightarrow 
\mathbb{C}$,
where we used the notations introduced in the proof of Proposition 
\ref{fusedAW}. The latter integral is essentially the integral
of $f$ against the Askey-Wilson measure 
depending on the four Askey-Wilson parameters $t_j$
($j=1,\ldots,4$) and with base $\widetilde{q}$. 
In particular, the Askey-Wilson polynomials
\cite{AW}, \cite[\S 7.5]{GR}
\[
p_n(x)={}_4\phi_3\left(\begin{matrix} \widetilde{q}\,{}^{-n}, 
\widetilde{q}\,{}^{n-1}t_1t_2t_3t_4, 
t_1\exp(2\pi ix),t_1\exp(-2\pi ix)\\
t_1t_2, t_1t_3, t_1t_4 \end{matrix};\,\widetilde{q},\widetilde{q}\,
\right),\quad n\in \mathbb{Z}_{\geq 0}
\]
form a basis of the space $\mathbb{C}\lbrack \cos(2\pi x)\rbrack$
of polynomials in $\cos(2\pi x)$, and they 
are orthogonal with respect to the complex pairing
\[\langle p_1,p_2\rangle=\int_{-\infty}^{\infty}
p_1(x)p_2(x)J(x)dx.
\]
The ``quadratic norms'' $\langle p_n,p_n\rangle$ 
($n\in \mathbb{Z}_{\geq 0}$) can be explicitly
evaluated using the quadratic norm evaluations of the Askey-Wilson 
polynomials, see \cite{AW} and \cite[\S 7.5]{GR} 
(the special case $n=0$ is the statement
of Proposition \ref{fusedAW}).
\end{rema}

The proposition is also valid
for $(s_1,s_2,s_3,s_4)=(q,-q,q^{\frac{1}{2}},0)$, in which case the resulting
integral
\begin{equation*}
\begin{split}
\int_{-\infty}^{\infty}
&\frac{\bigl(\exp(\pm 2\pi ix),-\exp(\pm 2\pi ix), 
\widetilde{q}\,{}^{\frac{1}{2}}\exp(\pm 2\pi ix);
\widetilde{q}\,\bigr)_{\infty}}
{\prod_{j=1}^4\bigl(t_j\exp(\pm 2\pi ix);\widetilde{q}\,\bigr)_{\infty}}\,
q^{\frac{x^2}{2}}\,dx\\
&\qquad\qquad\qquad\qquad\qquad\qquad\qquad\qquad=
\frac{2}{\sqrt{-i\tau}}\frac{\bigl(t_1t_2t_3t_4;\widetilde{q}\,\bigr)_{\infty}}
{\prod_{1\leq k<m\leq 4}\bigl(t_kt_m;\widetilde{q}\,\bigr)_{\infty}}
\end{split}
\end{equation*}
is the fusion of the Jacobi inversion formula \eqref{inversion} and
the trigonometric Askey-Wilson integral \eqref{betaEuler2}.
This integral may be viewed as a one-variable 
Macdonald-Mehta integral, see \cite{S2} and \cite{C}.
In particular, the special case
$(t_1,t_2,t_3,t_4)=\bigl(0,\widetilde{q}\,{}^{\frac{1}{2}},
\widetilde{q}\,{}^k, -\widetilde{q}\,{}^{k}\bigr)$
with $k\in\mathbb{R}_{>0}$ leads to Cherednik's \cite{C} one-variable
Macdonald-Mehta integral
\[\int_{-\infty}^{\infty}\frac{\bigl(\exp(\pm 4\pi ix);
\widetilde{q}\,{}^2\bigr)_{\infty}}
{\bigl(\widetilde{q}\,{}^{2k}\exp(\pm 4\pi ix);
\widetilde{q}\,{}^2\bigr)_{\infty}}\,q^{\frac{x^2}{2}}dx=
\frac{2}{\sqrt{-i\tau}}\frac{\bigl(\widetilde{q}\,{}^{2k};
\widetilde{q}\,{}^2\bigr)_{\infty}}
{\bigl(\widetilde{q}\,{}^{4k};\widetilde{q}\,{}^2\bigr)_{\infty}}.
\]
In a similar way we would like to put 
$(t_1,t_2,t_3,t_4)=(1,-1,\widetilde{q}\,{}^{\frac{1}{2}},0)$ in the
fused Askey-Wilson integral (Proposition \ref{fusedAW}), 
but this is not allowed since the
integrand is irregular on the real line for these parameter values.
A suitable regularization of this integral yields an
integral evaluation 
which is closely related to a trigonometric beta integral of
Askey \cite[\S 3]{A3}, and which
has a similar appearance as Etingof's \cite[\S 3.5]{CO}
Macdonald-Mehta
type integral with real integration cycle. The result is as follows.
\begin{cor}
Let $\tau\in\mathbb{H}$ and $s_j\in \mathbb{C}$ 
\textup{(}$j=1,\ldots,4$\textup{)} with $|q^{-3}s_1s_2s_3s_4|<1$.
Then 
\[
\lim_{\epsilon\searrow 0}\,
\int_{i\epsilon-\infty}^{i\epsilon+\infty}
\frac{\prod_{j=1}^4\bigl(s_jq^{\pm z};q\bigr)_{\infty}}
{\bigl(q^{1\pm z},-q^{1\pm z},q^{\frac{1}{2}\pm z};q\bigr)_{\infty}}
\,q^{\frac{z^2}{2}}\,dz
=\frac{1}{\sqrt{-i\tau}}
\frac{\prod_{1\leq k<m\leq 4}\bigl(q^{-1}s_ks_m;q\bigr)_{\infty}}
{\bigl(q^{-3}s_1s_2s_3s_4;q\bigr)_{\infty}}.
\]
\end{cor}
\begin{proof} 
The integrand $\widetilde{\phi}(z)$ is regular on $i\epsilon+\mathbb{R}$
for $\epsilon\in\mathbb{R}_{>0}$ sufficiently small.
Folding the integral yields 
\[
\int_{i\epsilon-\infty}^{i\epsilon+\infty}
\widetilde{\phi}(z)\,dz=
\int_0^1\phi^+(i\epsilon+x)dx=
\Phi\int_0^1\vartheta_{-\frac{1}{\tau}}(i\epsilon+x)dx
\]
by \eqref{Bailey} and \eqref{PhiAW}.
By \eqref{thetaintegral} we have
\[\lim_{\epsilon\searrow 0}\,\int_0^1\vartheta_{-\frac{1}{\tau}}(i\epsilon+x)dx
=\int_0^1\vartheta_{-\frac{1}{\tau}}(x)dx=1,
\]
which completes the proof.
\end{proof}

\subsection{The hyperbolic Askey-Wilson integral}

In this section we derive the hyperbolic analogue of the Askey-Wilson
integral and, as a special case, the analogue of the
Askey-Wilson integral for $|q|=1$. The techniques are the same as for 
the derivation of the hyperbolic Ramanujan integral in \S 3.4.

For the moment we assume that
$\tau\in\mathbb{H}\cap \mathbb{C}_-$.
As in the previous section we write
$q=q_\tau=\exp(2\pi i\tau)$ and $\widetilde{q}=\widetilde{q}_\tau=
\exp(-2\pi i/\tau)$, and $q^u=\exp(2\pi i\tau u)$, $\widetilde{q}\,{}^u=
\exp(-2\pi iu/\tau)$ for $u\in \mathbb{C}$.
We fix 
four parameters $\tau_j$ ($j=1,\ldots,4$) satisfying the parameter constraints
\begin{equation}\label{pardomainAWa}
\begin{split}
&1-\tau_j\in \mathbb{H}\cap \mathbb{C}_+,\quad
(1-\tau_j)\tau\in \mathbb{H}\quad (j=1,\ldots,4),\\
&\bigl(\tau_1+\tau_2+\tau_3+\tau_4-3-\tau^{-1}\bigr)\tau
\in\mathbb{H}\cap\mathbb{C}_-.
\end{split}
\end{equation}
If we take parameters $\tau_j\in\mathbb{C}$ satisfying 
\[\hbox{Re}(\tau_j)<1,\qquad \hbox{Im}(\tau_j)<0,\qquad
\hbox{Re}(\tau_1+\tau_2+\tau_3+\tau_4)>3,
\]
then sufficient additional conditions on their imaginary parts 
to satisfy the parameter constraints \eqref{pardomainAWa} are
\begin{equation}\label{nonempty}
\hbox{Im}(\tau_j)>\bigl(1-\hbox{Re}(\tau_j)\bigr)\frac{\hbox{Im}(\tau)}
{\hbox{Re}(\tau)},\qquad
\hbox{Im}(\tau_1+\tau_2+\tau_3+\tau_4)>\hbox{Re}\bigl(\tau_1+\tau_2+
\tau_3+\tau_4-3\bigr)\frac{\hbox{Re}(\tau)}{\hbox{Im}(\tau)}.
\end{equation}
The existence of parameters satisfying these inequalities follows
from the fact that the right hand sides of the inequalities in 
\eqref{nonempty} are strictly negative.

For fixed parameters $\tau_j$ satisfying the parameter constraints
\eqref{pardomainAWa} 
we associate eight parameters $s_j,t_j$ ($j=1,\ldots,4$) by
\begin{equation}\label{parspecializationAW}
s_j=q^{\tau_j},\qquad t_j=\exp(-2\pi i\tau_j),\qquad (j=1,\ldots,4).
\end{equation}
These eight parameters satisfy the parameter requirements of
Proposition \ref{fusedAW}. The resulting fused Askey-Wilson integral
can then be written in terms of the $\tau$-shifted factorial 
\eqref{taufactorial} as
\begin{equation}\label{step1AW}
\int_{-\infty}^{\infty}\frac{\lbrack 1\pm x, \frac{1}{2}\pm x,
1+\frac{1}{2\tau}\pm x;\tau\rbrack_{\infty}}
{\prod_{j=1}^4\lbrack \tau_j\pm x;\tau\rbrack_{\infty}}\,q^{\frac{x^2}{2}}\,dx
=\frac{2}{\sqrt{-i\tau}}\frac{\lbrack \tau_1+\tau_2+\tau_3+\tau_4-3;
\tau\rbrack_{\infty}}
{\prod_{1\leq k<m \leq 4}\lbrack \tau_k+\tau_m-1;\tau\rbrack_{\infty}}.
\end{equation}
Using the reflection equation \eqref{reflectioneq} and 
\begin{equation}\label{doubling}
\lbrack \frac{1}{2}\pm x, 1\pm x, \frac{1}{2}+\frac{1}{2\tau}\pm x,
1+\frac{1}{2\tau}\pm x;\tau\rbrack_{\infty}=
\lbrack 1\pm 2x;\tau\rbrack_{\infty}
\end{equation}
(which follows easily from \eqref{taufactorial} and \eqref{expanding}),
we can rewrite \eqref{step1AW} as
\begin{equation}\label{step2AW}
\int_{-\infty}^\infty
\frac{\lbrack 1\pm 2x;\tau\rbrack_{\infty}}
{\prod_{j=1}^4\lbrack \tau_j\pm x;\tau\rbrack_{\infty}}\,dx\\
=\frac{2q^{-\frac{1}{24}}\widetilde{q}\,{}^{\frac{1}{24}}}{\sqrt{-i\tau}}\,
\frac{\lbrack \tau_1+\tau_2+\tau_3+\tau_4-3;\tau\rbrack_{\infty}}
{\prod_{1\leq k<m\leq 4}\lbrack \tau_k+\tau_m-1;\tau\rbrack_{\infty}},
\end{equation}
which now has a very similar appearance as 
the trigonometric Askey-Wilson integral \eqref{betaEuler2}.
The numerator of the integrand
\begin{equation}\label{IAW}
I(z)=\frac{\lbrack 1\pm 2z;\tau\rbrack_{\infty}}
{\prod_{j=1}^4\lbrack \tau_j\pm z;\tau\rbrack_{\infty}}
\end{equation}
can be simplified as follows.
\begin{lem}\label{numerator}
\[\lbrack 1\pm 2z;\tau\rbrack_{\infty}=
-iq^{\frac{1}{12}}\widetilde{q}\,{}^{-\frac{1}{12}}q^{2z^2}
\bigl(q^{-z}-q^z\bigr)\bigl(\exp(2\pi iz)-\exp(-2\pi iz)\bigr).
\]
\end{lem}
\begin{proof}
By the reflection equation \eqref{reflectioneq} we can write
\[\lbrack 1\pm 2z;\tau\rbrack_{\infty}=
q^{-\frac{1}{24}}\widetilde{q}\,{}^{\frac{1}{24}}
q^{\frac{1}{2}\bigl(\frac{1}{2}-\frac{1}{2\tau}-2z\bigr)^2}\,
\frac{\lbrack 2z+1;\tau\rbrack_{\infty}}
{\lbrack 2z+\tau^{-1};\tau\rbrack_{\infty}}.
\]
Now apply the functional equations \eqref{functionalequations}
to the numerator and denominator to get rid of the $\tau$-shifted
factorials. Simplifying the resulting elementary expression gives
the desired result.
\end{proof}
By the zero and pole locations of the $\tau$-shifted factorial
(see \eqref{zero} and \eqref{pole} respectively) we conclude that
the poles of $I(z)$ are contained in the eight discrete sets
\begin{equation}\label{poleIAW}
\pm\left(\tau_j+\frac{1}{\tau}\mathbb{Z}_{\geq 0}+\mathbb{Z}_{<0}\right)
\qquad (j=1,\ldots,4).
\end{equation}
By the parameter constraints $\tau\in\mathbb{H}\cap \mathbb{C}_-$
and \eqref{pardomainAWa}, the four sets in \eqref{poleIAW}
with plus sign are contained in the interior of the third quadrant of
the complex plane, and consequently the four sets with minus sign are
contained in the interior of the first quadrant of the complex plane.

\begin{lem}\label{boundAW}
Suppose that $\tau\in \mathbb{H}\cap \mathbb{C}_-$ and that the parameters
$\tau_j$ \textup{(}$j=1,\ldots,4$\textup{)} satisfy \eqref{pardomainAWa}.
Then there exists a constant $C\in \mathbb{R}_{>0}$ such that
\[|I(z)|\leq C\left|
q^{(\tau_1+\tau_2+\tau_3+\tau_4-3-\tau^{-1})z}
\right|
\]
for $z\in \mathbb{C}$ satisfying
$\hbox{Re}(z)\geq 0$ and $\hbox{Im}(z)\leq 0$.
\end{lem}
\begin{proof}
By Lemma \ref{numerator}, the reflection equation \eqref{reflectioneq}
and by \eqref{taufactorial} we can rewrite the integrand $I(z)$ as
\begin{equation*}
\begin{split}
I(z)=&C
\prod_{j=1}^4\frac{\bigl(\widetilde{q}\,\exp(2\pi i\tau_j)\exp(-2\pi iz);
\widetilde{q}\,\bigr)_{\infty}}
{\bigl(\exp(-2\pi i\tau_j)\exp(-2\pi iz);\widetilde{q}\,\bigr)_{\infty}}
\frac{\bigl(q^{\tau_j}q^z;q\bigr)_{\infty}}
{\bigl(q^{1-\tau_j}q^z;q\bigr)_{\infty}}\\
&\qquad\qquad\times \bigl(1-\exp(-4\pi iz)\bigr)\bigl(1-q^{2z}\bigr)
\,q^{(\tau_1+\tau_2+\tau_3+\tau_4-3-\tau^{-1})z}
\end{split}
\end{equation*}
for some irrelevant nonzero constant $C$. The proof is now analogous to the
proof of Lemma \ref{boundRam}, since $|\exp(-2\pi i\tau_j)|<1$
and $|q^{1-\tau_j}|<1$ by the conditions \eqref{pardomainAWa} on
the parameters. 
\end{proof}
Note that a bound for the weight function $I(z)$ for $z\in \mathbb{C}$
satisfying $\hbox{Re}(z)\leq 0$
and $\hbox{Im}(z)\geq 0$ can be immediately deduced from Lemma \ref{boundAW}
using $I(-z)=I(z)$. The following corollary follows from
Lemma \ref{boundAW} and the conditions \eqref{pardomainAWa} on the parameters.
\begin{cor}\label{corboundAW}
Let $\tau\in \mathbb{H}\cap \mathbb{C}_-$ and suppose that the parameters
$\tau_j$ satisfy the parameter constraints \eqref{pardomainAWa}.
Let $0<\epsilon<1$ be the growth coefficient
$\epsilon=\exp(-2\pi m)$ with $m$ the strictly positive constant
\[
m=\min_{\theta\in \lbrack -\frac{\pi}{2},0\rbrack}
\left(\hbox{Im}\bigl((\tau_1+\tau_2+\tau_3+\tau_4-3-\tau^{-1})
\tau e^{i\theta}\bigr)
\right).
\]
Then there exists a constant $C\in \mathbb{R}_{>0}$ such that
\[|I(z)|\leq C\epsilon^{|z|}
\]
for $z\in \mathbb{C}$ satisfying  $\hbox{Re}(z)\geq 0$ and 
$\hbox{Im}(z)\leq 0$, as well as for $z\in\mathbb{C}$ satisfying
$\hbox{Re}(z)\leq 0$ and $\hbox{Im}(z)\geq 0$.
\end{cor}
We keep the assumptions $\tau\in \mathbb{H}\cap \mathbb{C}_-$
and \eqref{pardomainAWa}. 
In view of Corollary \ref{corboundAW} we may apply
Cauchy's theorem to rotate clockwise the integration
contour $\mathbb{R}$ in \eqref{step2AW} to $-i\mathbb{R}$ about the origin
without altering its evaluation.
Relaxing the parameter constraints
leads now to the following result.
\begin{thm}[Hyperbolic Askey-Wilson integral]\label{hyperbolicAW}
Let $\tau\in \mathbb{C}$ with $\hbox{Re}(\tau)<0$ and 
$\hbox{Im}(\tau)\geq 0$, and suppose that the parameters
$\tau_j\in \mathbb{C}$ \textup{(}$j=1,\ldots,4$\textup{)} satisfy
\begin{equation}\label{pardomainAWb}
\hbox{Re}(\tau_j)<1\quad (j=1,\ldots,4),\qquad
\hbox{Re}\bigl((\tau_1+\tau_2+\tau_3+\tau_4-3)\tau\bigr)<1.
\end{equation}
Then
\begin{equation}\label{AW3}
\int_{-i\infty}^{i\infty}
\frac{\lbrack 1\pm 2z;\tau\rbrack_{\infty}}
{\prod_{j=1}^4\lbrack \tau_j\pm z;\tau\rbrack_{\infty}}\,dz
=-\frac{2q^{-\frac{1}{24}}\widetilde{q}\,{}^{\frac{1}{24}}}{\sqrt{-i\tau}}\,
\frac{\lbrack \tau_1+\tau_2+\tau_3+\tau_4-3;\tau\rbrack_{\infty}}
{\prod_{1\leq k<m\leq 4}\lbrack \tau_k+\tau_m-1;\tau\rbrack_{\infty}}.
\end{equation}
\end{thm}
\begin{proof}
For $\tau\in \mathbb{H}\cap \mathbb{C}_-$ we already argued that
the integral evaluation \eqref{AW3} is valid for parameters $\tau_j$
satisfying the parameter constraints \eqref{pardomainAWa}. 
By analytic continuation one proves that the integral evaluation \eqref{AW3} 
is valid for parameters $\tau_j$ satisfying the milder constraints
\eqref{pardomainAWb}. In fact, by \eqref{as1} and \eqref{as2} 
the open convex parameter domain of parameters 
$(\tau_1,\tau_2,\tau_3,\tau_4)$ satisfying the constraints 
\eqref{pardomainAWb} is the maximal extension of the parameter 
domain \eqref{pardomainAWa}
such that the associated integrand $I(z)$ decays exponentially as 
$z\rightarrow\pm i\infty$ and such that the
sequence of poles $\tau_j+\frac{1}{\tau}\mathbb{Z}_{\geq 0}+\mathbb{Z}_{<0}$
($j=1,\ldots,4$) lie in $\mathbb{C}_-$.

The extension of the integral evaluation \eqref{AW3}
to $\tau\in \mathbb{R}_{<0}$ follows now in a similar manner
as for the hyperbolic Ramanujan integral (see the proof of 
Theorem \ref{hyperbolicR}) using the  
reflection equation \eqref{reflectioneq} and the   
asymptotics \eqref{as1} and \eqref{as2}
of the $\tau$-shifted factorial.
\end{proof}
\begin{rema}\label{AWpol}
The Askey-Wilson polynomials (see Remark \ref{polAWremark})
do not remain an orthogonal system in the hyperbolic case. 
Indeed, the Askey-Wilson polynomials are one-periodic hence 
bounded on $\mathbb{R}$, but they 
grow exponentially in the imaginary direction (with growth rate
$\exp(2\pi n|\hbox{Im}(z)|)$ where 
$n$ is the degree of the Askey-Wilson polynomial). 
The rotation of the integration cycle $\mathbb{R}$ to $-i\mathbb{R}$
in the complex orthogonality relations of Remark \ref{polAWremark}
is only allowed if the integrand has uniform 
exponential decay in the second and fourth quadrant of the complex plane. 
For given parameters satisfying \eqref{pardomainAWa},
this will thus only be the case for certain 
low degree Askey-Wilson polynomials. 

On the other hand, in \cite{S} an hyperbolic 
analogue of the Askey-Wilson function is defined, which is expected to
play the role of integral kernel for some generalized $q$-Fourier transform
with $|q|=1$.
\end{rema}

\section{The hyperbolic Nassrallah-Rahman integral}


\subsection{A weak ${}_8\Psi_8$ summation formula}

We take $\tau\in\mathbb{H}$, so that $|q|<1$, and we fix five generic
nonzero complex parameters $t_j$ ($j=0,\ldots,4$). It is convenient to
write $A=t_0t_1t_2t_3t_4$ for the product of the five parameters.
For the hyperbolic analogue of the 
Nassrallah-Rahman integral, the role
of the function $\phi$ in the underlying summation formula (see \S 2)
is 
\begin{equation}\label{phi3}
\phi(z)=
\frac{\prod_{j=0}^4\bigl(t_jq^{\pm z};q\bigr)_{\infty}}
{\bigl(q^{1\pm z}, -q^{1\pm z},   
q^{\frac{1}{2}\pm z},Aq^{-4\pm z};q\bigr)_{\infty}}.
\end{equation}
As in \S 2, we define 
\[
\widetilde{\phi}(z)=\phi(z)\,q^{\frac{z^2}{2}}.
\]
A direct computation shows that
\[\widetilde{\phi}(z+1)=t(z)\widetilde{\phi}(z)
\]
with $t(z)$ given by
\[t(z)=q\frac{\bigl(1-q^{1+z}\bigr)\bigl(1+q^{1+z}\bigr)
\bigl(1-Aq^{-4+z}\bigr)}
{\bigl(1-q^{z}\bigr)\bigl(1+q^{z}\bigr)\bigl(1-
A^{-1}q^{5+z}\bigr)}
\prod_{j=0}^4\frac{\bigl(1-t_j^{-1}q^{1+z}\bigr)}
{\bigl(1-t_jq^z\bigr)}.
\]
The bilateral sum 
\begin{equation*}
\begin{split}
t^+(z)=&\sum_{m=0}^{\infty}\prod_{k=0}^{m-1}t(z+k)+
\sum_{m=-\infty}^{-1}\prod_{k=m}^{-1}\frac{1}{t(z+k)}\\
=&{}_8\Psi_8\left(\begin{matrix} q^{1+z}, -q^{1+z}, t_0^{-1}q^{1+z}, 
t_1^{-1}q^{1+z},
t_2^{-1}q^{1+z}, t_3^{-1}q^{1+z}, t_4^{-1}q^{1+z}, 
Aq^{-4+z}\\
q^z, -q^z, t_0q^z, t_1q^z, t_2q^z, t_3q^z, t_4q^z,
A^{-1}q^{5+z}
\end{matrix};q,q\right)
\end{split}
\end{equation*}
converges absolutely and uniformly on compacta away from the 
$\mathbb{Z}_{\leq 0}$-translates of the poles and the 
$\mathbb{Z}_{>0}$-translates of the zeros of $t(z)$, and
it extends to a meromorphic function in $z\in\mathbb{C}$.
Hence the bilateral sum
\[\phi^+(z)=\sum_{m=-\infty}^{\infty}\widetilde{\phi}(z+m)
\]
converges to a meromorphic function in $z\in\mathbb{C}$, given explicitly 
by
\[\phi^+(z)=t^+(z)\widetilde{\phi}(z).
\]
Unfortunately, there is no explicit summation formula
on the ${}_8\Psi_8$ level. Instead we apply a formula
that expresses $t^+(z)$ as a combination of two one-sided sums 
which only involve $z$-independent summands. 
This can be done by applying \cite[(5.6.2)]{GR} with its seven parameters
$a,b,c,\ldots,g$ specialized to
\[(a,b,c,d,e,f,g)=
\bigl(q^z,qt_1^{-1},qt_2^{-1},qt_3^{-1},qt_4^{-1},qt_0^{-1},
q^{-4}A\bigr).
\]
After a straightforward but tedious
computation using \eqref{expanding}, the Jacobi triple product
identity and the Jacobi inversion formula, we then arrive at the explicit
expression
\[
t^+(z)=\left\{C_1+\frac{\bigl(t_0q^{\pm z},t_0^{-1}q^{1\pm z};
q\bigr)_{\infty}}
{\bigl(Aq^{-4\pm z},A^{-1}q^{5\pm z};q\bigr)_{\infty}}C_2\right\}
\widetilde{\phi}(z)^{-1}\vartheta_{-\frac{1}{\tau}}(z)
\]
with the $z$-independent constants $C_1$ and $C_2$ given by
\begin{equation}\label{C1}
\begin{split}
C_1=\frac{1}{\sqrt{-i\tau}}&\frac{\prod_{j=1}^4\bigl(q^{-1}t_0t_j,
qt_0^{-1}t_j;q\bigr)_{\infty}}
{\bigl(q^{-5}At_0,q^{-3}At_0^{-1},
q^3t_0^{-2};q\bigr)_{\infty}}\\
&\times
{}_8W_7\bigl(q^2t_0^{-2};q^2t_0^{-1}t_1^{-1}, q^2t_0^{-1}t_2^{-1},
q^2t_0^{-1}t_3^{-1}, q^2t_0^{-1}t_4^{-1},q^{-3}At_0^{-1};q,q\bigr) 
\end{split}
\end{equation}
and
\begin{equation}\label{C2}
\begin{split}
C_2=\frac{1}{\sqrt{-i\tau}}&\frac{\prod_{j=1}^4
\bigl(q^4A^{-1}t_j,q^{-4}At_j;q\bigr)_{\infty}}
{\bigl(q^5A^{-1}t_0^{-1},q^{-3}At_0^{-1},q^{-7}A^2;q\bigr)_{\infty}}\\
&\times{}_8W_7\bigl(q^{-8}A^2;q^{-3}At_0^{-1},
q^{-3}At_1^{-1},q^{-3}At_2^{-1},q^{-3}At_3^{-1},q^{-3}At_4^{-1};q,q\bigr).
\end{split}
\end{equation}
The meromorphic function $\phi^+(z)$ can thus be written as
\begin{equation}\label{8Psi8}
\phi^+(z)=\Phi(z)\vartheta_{-\frac{1}{\tau}}(z),\qquad
\Phi(z)=C_1+\frac{\bigl(t_0q^{\pm z},t_0^{-1}q^{1\pm z};
q\bigr)_{\infty}}
{\bigl(Aq^{-4\pm z},A^{-1}q^{5\pm z};q\bigr)_{\infty}}\,C_2.
\end{equation}
The function $\Phi(z)$ is an elliptic function in $z$ with respect to the
periods $1$ and $-\frac{1}{\tau}$.
For later purposes it is convenient to rewrite $\Phi(z)$ as follows.
Fix $\tau_j\in \mathbb{C}$ ($j=0,\ldots,4$) such that
\[
t_j=q^{\tau_j}\qquad (j=0,\ldots,4)
\]
and denote $a=\tau_0+\tau_1+\tau_2+\tau_3+\tau_4$, so that
$A=q^a$. We define new nonzero complex parameters $\widetilde{t}_j$
($j=0,\ldots,4$) by
\[
\widetilde{t}_j=\exp(-2\pi i\tau_j)\qquad (j=0,\ldots,4),
\]
and we write $\widetilde{A}=
\widetilde{t}_0\widetilde{t}_1\widetilde{t}_2\widetilde{t}_3\widetilde{t}_4=
\exp(-2\pi ia)$.
Then
\begin{equation}\label{PhiNR}
\begin{split}
\Phi(z)&=C_1+\,\frac{\bigl(\widetilde{t}_0\exp(\pm 2\pi iz),
\widetilde{q}\,\widetilde{t}_0^{-1}\exp(\pm 2\pi iz);\widetilde{q}\,
\bigr)_{\infty}}
{\bigl(\widetilde{A}\,\exp(\pm 2\pi iz), \widetilde{q}\,
\widetilde{A}^{-1}\exp(\pm 2\pi iz);\widetilde{q}\,\bigr)_{\infty}}\, 
C_2q^\gamma,\\
\gamma&=\left(a-\frac{9}{2}-\frac{1}{2\tau}\right)^2-
\left(\tau_0-\frac{1}{2}-\frac{1}{2\tau}\right)^2
\end{split}
\end{equation}
by a direct computation using the Jacobi triple product identity
and the Jacobi inversion formula.


\subsection{Fusion} In this subsection we fuse the weak ${}_8\Psi_8$
summation formula \eqref{8Psi8} with the Nassrallah-Rahman integral
\eqref{betaEuler2n}. For proofs of the Nassrallah-Rahman integral 
\eqref{betaEuler2n} itself, see \cite{NR}, \cite[(6.4.1)]{GR} 
and references therein.

We fix $\tau\in \mathbb{H}$ and generic complex parameters $\tau_j$. 
As in the previous subsection, we write 
\begin{equation}\label{fusedparametersNR}
t_j=q^{\tau_j},\qquad \widetilde{t}_j=\exp(-2\pi i\tau_j)\qquad 
\qquad (j=0,\ldots,4)
\end{equation}
and $A=q^a$, $\widetilde{A}=\exp(-2\pi ia)$ with 
$a=\tau_0+\tau_1+\tau_2+\tau_3+\tau_4$. We assume throughout this subsection
that  $\hbox{Im}(\tau_j)<0$ for $j=0,\ldots,4$ and that
$\hbox{Im}(a)>\hbox{Im}(\tau^{-1})$. 
These conditions imply that 
\[|\widetilde{t}_j|<1\quad (j=0,\ldots,4),\qquad\qquad
|\widetilde{A}\,|>|\widetilde{q}\,|.
\]
By \eqref{doubling}, 
the reflection equation \eqref{reflectioneq}
and the definition \eqref{taufactorial}
of the $\tau$-shifted factorial, we can write 
\begin{equation}\label{rewritten}
\frac{\lbrack 1\pm 2z, a-4\pm z;\tau\rbrack_{\infty}}
{\prod_{j=0}^4\lbrack \tau_j\pm z;\tau\rbrack_{\infty}}=
q^{-\frac{1}{24}}\widetilde{q}\,{}^{\frac{1}{24}}\widetilde{\phi}(z)
v(z)
\end{equation}
with $\widetilde{\phi}(z)$
as in \S 5.1, and with $v(z)$ the one-periodic
function
\begin{equation}
v(z)=\frac{\bigl(\exp(\pm 2\pi iz),-\exp(\pm 2\pi iz),
\widetilde{q}\,{}^{\frac{1}{2}}\exp(\pm 2\pi iz),
\widetilde{A}\,\exp(\pm 2\pi iz);\widetilde{q}\,\bigr)_{\infty}}
{\prod_{j=0}^4
\bigl(\widetilde{t}_j\exp(\pm 2\pi iz);\widetilde{q}\,\bigr)_{\infty}}.
\end{equation}
Note that \eqref{rewritten} is regular on $\mathbb{R}$ in view
of the conditions on the parameters $\tau_j$.
By Fubini's theorem, the integral
\begin{equation}
\int_{-\infty}^{\infty}\frac{\lbrack 1\pm 2x, a-4\pm x;\tau\rbrack_{\infty}}
{\prod_{j=0}^4\lbrack \tau_j\pm x;\tau\rbrack_{\infty}}\,dx=
q^{-\frac{1}{24}}\widetilde{q}\,{}^{\frac{1}{24}}
\int_{-\infty}^{\infty}\widetilde{\phi}(x)v(x)\,dx
\end{equation}
can now be folded,
\begin{equation*}
\begin{split}
\int_{-\infty}^{\infty}\widetilde{\phi}(x)v(x)\,dx
&=\int_0^1\phi^+(x)v(x)\,dx\\
&=C_1\int_0^1v(x)\vartheta_{-\frac{1}{\tau}}(x)\,dx\\
&+C_2q^\gamma\,
\int_0^1v(x)\frac{\bigl(\widetilde{t}_0\exp(\pm 2\pi ix),
\widetilde{q}\,
\widetilde{t}_0^{-1}\exp(\pm 2\pi ix);\widetilde{q}\,\bigr)_{\infty}}
{\bigl(\widetilde{A}\exp(\pm 2\pi ix), 
\widetilde{q}\,
\widetilde{A}^{-1}\exp(\pm 2\pi ix);\widetilde{q}\,\bigr)_{\infty}}
\vartheta_{-\frac{1}{\tau}}(x)\,dx
\end{split}
\end{equation*}
by \eqref{8Psi8} and \eqref{PhiNR}.
Applying the Jacobi triple product identity and \eqref{expanding}, 
the remaining two integrals over the period interval $[0,1]$ 
can be evaluated by the trigonometric
Nassrallah-Rahman integral \eqref{betaEuler2n}. This gives
\begin{equation}\label{tussenstapa}
\int_{-\infty}^{\infty}\widetilde{\phi}(x)v(x)\,dx=
\frac{2\prod_{j=0}^4\bigl(\widetilde{A}\,
\widetilde{t}_j^{-1};\widetilde{q}\,\bigr)_{\infty}}
{\prod_{0\leq k<m\leq 4} 
\bigl(\widetilde{t}_k\widetilde{t}_m;\widetilde{q}\,\bigr)_{\infty}}
\left\{C_1+C_2q^{\gamma}\,
\prod_{j=1}^4
\frac{\bigl(\widetilde{t}_0\widetilde{t}_j,
\widetilde{q}\,\widetilde{t}_0^{-1}\widetilde{t}_j^{-1};
\widetilde{q}\,\bigr)_{\infty}}
{\bigl(\widetilde{A}\,\widetilde{t}_j^{-1},
\widetilde{q}\,\widetilde{A}^{-1}\widetilde{t}_j;
\widetilde{q}\,\bigr)_{\infty}}\right\},
\end{equation}
with the constants $C_1$ and $C_2$ given by \eqref{C1} and \eqref{C2}, 
respectively. The following key lemma gives the evaluation of  
the remaining sum.
\begin{lem}
\[C_1+C_2q^\gamma\,\prod_{j=1}^4\frac{\bigl(\widetilde{t}_0\widetilde{t}_j,
\widetilde{q}\,\widetilde{t}_0^{-1}\widetilde{t}_j^{-1};
\widetilde{q}\,\bigr)_{\infty}}
{\bigl(\widetilde{A}\,\widetilde{t}_j^{-1},
\widetilde{q}\,\widetilde{A}^{-1}\widetilde{t}_j;
\widetilde{q}\,\bigr)_{\infty}}=\frac{1}{\sqrt{-i\tau}}
\frac{\prod_{0\leq k<m\leq 4}\bigl(q^{-1}t_kt_m;q\bigr)_{\infty}}
{\prod_{j=0}^4\bigl(q^{-3}At_j^{-1};q\bigr)_{\infty}}.
\]
\end{lem}
\begin{proof}
Using the Jacobi triple product identity and the
Jacobi inversion formula we can write
\[\prod_{j=1}^4\frac{\bigl(\widetilde{t}_0\widetilde{t}_j,
\widetilde{q}\,\widetilde{t}_0^{-1}\widetilde{t}_j^{-1};
\widetilde{q}\,\bigr)_{\infty}}
{\bigl(\widetilde{A}\,\widetilde{t}_j^{-1},
\widetilde{q}\,\widetilde{A}^{-1}\widetilde{t}_j;
\widetilde{q}\,\bigr)_{\infty}}=
q^\delta\,\prod_{j=1}^4
\frac{\bigl(q^{-1}t_0t_j,q^2t_0^{-1}t_j^{-1};q\bigr)_{\infty}}
{\bigl(q^{-3}At_j^{-1}, q^4A^{-1}t_j;q\bigr)_{\infty}}
\]
with $\delta$ given explicitly by
\[\delta=\frac{1}{2}\sum_{j=1}^4\left\{
\left(\tau_0+\tau_j-\frac{3}{2}-\frac{1}{2\tau}\right)^2-
\left(a-\tau_j-\frac{7}{2}-\frac{1}{2\tau}\right)^2\right\}.
\]
Writing out the squares for $\gamma$ \eqref{PhiNR} and $\delta$ gives
\[\gamma=4\tau^{-1}+20-8\tau_0+
\sum_{j=1}^4\bigl(\tau_j^2-(9+\tau^{-1})\tau_j\bigr)
+2\sum_{0\leq k<m\leq 4}\tau_k\tau_m=-\delta,
\]
hence 
\begin{equation}\label{hhhhh}
C_1+C_2q^\gamma\,\prod_{j=1}^4\frac{\bigl(\widetilde{t}_0\widetilde{t}_j,
\widetilde{q}\,\widetilde{t}_0^{-1}\widetilde{t}_j^{-1};
\widetilde{q}\,\bigr)_{\infty}}
{\bigl(\widetilde{A}\,\widetilde{t}_j^{-1},
\widetilde{q}\,\widetilde{A}^{-1}\widetilde{t}_j;
\widetilde{q}\,\bigr)_{\infty}}=
C_1+C_2\prod_{j=1}^4
\frac{\bigl(q^{-1}t_0t_j,q^2t_0^{-1}t_j^{-1};q\bigr)_{\infty}}
{\bigl(q^{-3}At_j^{-1}, q^4A^{-1}t_j;q\bigr)_{\infty}}
\end{equation}
with $C_1$ and $C_2$ given by \eqref{C1} and \eqref{C2}, respectively.
The right hand side of \eqref{hhhhh}
can be evaluated by Bailey's summation formula \cite[(2.11.7)]{GR}
with the six parameters $a,b,\ldots,f$ specialized to
\[(a,b,c,d,e,f)=(q^2t_0^{-2}, q^{-3}At_0^{-1},
q^2t_0^{-1}t_1^{-1},q^2t_0^{-1}t_2^{-1},q^2t_0^{-1}t_3^{-1},
q^2t_0^{-1}t_4^{-1}).
\] 
This yields the desired result.
\end{proof}
Combining the lemma with \eqref{rewritten}
and \eqref{tussenstapa} immediately implies the following result.
\begin{prop}\label{fusionNR}
Let $\tau\in \mathbb{H}$ and $\tau_j\in \mathbb{C}$ with
$\hbox{Im}(\tau_j)<0$ \textup{(}$j=0,\ldots,4$\textup{)}
and $\hbox{Im}(a)>\hbox{Im}(\tau^{-1})$, 
where $a=\tau_0+\tau_1+\tau_2+\tau_3+\tau_4$.
Then
\begin{equation}\label{NRfusedformula}
\int_{-\infty}^{\infty}
\frac{\lbrack 1\pm 2x,a-4\pm x;\tau\rbrack_{\infty}}
{\prod_{j=0}^4\lbrack \tau_j\pm x;\tau\rbrack_{\infty}}\,dx=
\frac{2q^{-\frac{1}{24}}
\widetilde{q}\,{}^{\frac{1}{24}}}{\sqrt{-i\tau}}\,
\frac{\prod_{j=0}^4\lbrack a-\tau_j-3;\tau\rbrack_{\infty}}
{\prod_{0\leq k<m\leq 4}\lbrack \tau_k+\tau_m-1;\tau\rbrack_{\infty}}.
\end{equation}
\end{prop}
\begin{rema}
The weak ${}_8\Psi_8$ summation formula \eqref{8Psi8}
involves an elliptic function $\Phi(z)$ depending
nontrivially on the parameters $t_j$. Consequently, the fusion of
the weak ${}_8\Psi_8$ summation formula with
the trigonometric Nassrallah-Rahman integral \eqref{betaEuler2n}
only has an explicit evaluation when the parameters in the 
trigonometric Nassrallah-Rahman integral are matched
to the parameters in the weak ${}_8\Psi_8$ summation formula
via \eqref{fusedparametersNR}. This is an essential difference
with the fused Ramanujan integral in \S 3.3 and the fused Askey-Wilson
integral in \S 4.2.
\end{rema}
It is possible to fuse the trigonometric Nassrallah-Rahman integral
\eqref{betaEuler2n} in base $\widetilde{q}$
with the Jacobi inversion formula \eqref{inversion},
yielding the one-variable Macdonald-Mehta type integral
\begin{equation}
\begin{split}
\int_{-\infty}^{\infty}&\frac{\bigl(\exp(\pm 2\pi ix),
-\exp(\pm 2\pi ix), \widetilde{q}\,{}^{\frac{1}{2}}\exp(\pm 2\pi ix),
A\exp(\pm 2\pi ix);\widetilde{q}\,\bigr)_{\infty}}
{\prod_{j=0}^4\bigl(t_j\exp(\pm 2\pi ix);\widetilde{q}\,\bigr)_{\infty}}
\,q^{\frac{x^2}{2}}\,dx\\
&\qquad\qquad\qquad\qquad\qquad\qquad\qquad
=\frac{2}{\sqrt{-i\tau}}\,\frac{\prod_{j=0}^4
\bigl(At_j^{-1};\widetilde{q}\,\bigr)_{\infty}}
{\prod_{0\leq k<m\leq 4}\bigl(t_kt_m;\widetilde{q}\,\bigr)_{\infty}}
\end{split}
\end{equation}
 with $A=t_0t_1t_2t_3t_4$ and with parameters $t_j\in \mathbb{C}$
satisfying $|t_j|<1$ ($j=0,\ldots,4$).

\subsection{The hyperbolic Nassrallah-Rahman integral}
Let $\tau\in \mathbb{H}\cap \mathbb{C}_-$ with $\hbox{Re}(\tau^{-1})<-1$ 
and fix five parameters $\tau_j$ ($j=0,\ldots,4$) satisfying
\begin{equation}\label{pardomainNRa}
1-\tau_j,a-4-\tau^{-1}\in \mathbb{H}\cap \mathbb{C}_+,\qquad
(1-\tau_j)\tau, (a-4)\tau\in \mathbb{H}
\end{equation}
for $j=0,\ldots,4$, where  
$a=\tau_0+\tau_1+\tau_2+\tau_3+\tau_4$. If the real parts of the
parameters $\tau_j$ satisfy
\begin{equation}\label{Re}
\hbox{Re}(\tau_j)<1 \quad
(j=0,\ldots,4),\qquad \hbox{Re}(a)>4,
\end{equation}
then they satisfy \eqref{pardomainNRa} 
when the imaginary parts of $\tau_j$ satisfy 
\begin{equation}\label{Im}
\bigl(1-\hbox{Re}(\tau_j)\bigr)\frac{\hbox{Im}(\tau)}{\hbox{Re}(\tau)}<
\hbox{Im}(\tau_j)<0,\qquad
\hbox{Im}(a)\frac{\hbox{Re}(\tau)}{\hbox{Im}(\tau)}<-1-\hbox{Re}(\tau^{-1})
\end{equation}
for $j=0,\ldots,4$. Note that the condition $\hbox{Re}(\tau^{-1})<-1$
implies $-1-\hbox{Re}(\tau^{-1})>0$, hence there exist parameters $\tau_j$
satisfying \eqref{Re} and \eqref{Im}. 

The constraints \eqref{pardomainNRa} imply the parameter
conditions of Proposition \ref{fusionNR}. 
We now rotate the integration cycle $\mathbb{R}$ 
clockwise to $-i\mathbb{R}$ about the origin in the fused Nassrallah-Rahman
integral \eqref{NRfusedformula}.
By the zero and pole locations of the $\tau$-shifted factorial
(see \eqref{zero} and \eqref{pole} respectively) and by Lemma \ref{numerator},
the poles of the integrand
\begin{equation}\label{INR}
I(z)=\frac{\lbrack 1\pm 2z, a-4\pm z;\tau\rbrack_{\infty}}
{\prod_{j=0}^4\lbrack \tau_j\pm z;\tau\rbrack_{\infty}}
\end{equation}
of the fused Nassrallah-Rahman integral \eqref{NRfusedformula}
are contained in the union of the twelve discrete sets
\begin{equation}\label{poleINR}
\pm\left(\tau_j+\frac{1}{\tau}\mathbb{Z}_{\geq 0}+\mathbb{Z}_{<0}\right)\qquad
(j=0,\ldots,4),\qquad
\pm\left(-a+\frac{1}{\tau}\mathbb{Z}_{>0}+\mathbb{Z}_{\leq 4}\right).
\end{equation}
By the parameter constraints $\tau\in \mathbb{H}\cap \mathbb{C}_-$, 
$\hbox{Re}(\tau^{-1})<-1$ and
\eqref{pardomainNRa}, the six sets in \eqref{poleINR}
with plus sign are contained in the interior of the third quadrant
of the complex plane, and consequently the six sets in \eqref{poleINR}
with minus sign are contained in the interior of the first quadrant of
the complex plane. Furthermore, the integrand $I(z)$ is even, $I(-z)=I(z)$,
and by \eqref{reflectioneq} and 
\eqref{taufactorial} it can be rewritten as
\begin{equation*}
\begin{split}
I(z)=C\,&\bigl(1-q^{2z}\bigr)\bigl(1-\exp(-4\pi iz)\bigr)
\frac{\bigl(\exp(-2\pi ia)\exp(-2\pi iz);\widetilde{q}\,\bigr)_{\infty}}
{\bigl(\widetilde{q}\,\exp(2\pi ia)\exp(-2\pi iz);\widetilde{q}\,
\bigr)_{\infty}}
\frac{\bigl(q^{5-a}q^z;q\bigr)_{\infty}}
{\bigl(q^{a-4}q^z;q\bigr)_{\infty}}\\
&\quad\quad\times
\left\{\prod_{j=0}^4\frac{\bigl(\widetilde{q}\,
\exp(2\pi i\tau_j)\exp(-2\pi iz);\widetilde{q}\,\bigr)_{\infty}}
{\bigl(\exp(-2\pi i\tau_j)\exp(-2\pi iz);\widetilde{q}\,\bigr)_{\infty}}
\frac{\bigl(q^{\tau_j}q^z;q\bigr)_{\infty}}
{\bigl(q^{1-\tau_j}q^z;q\bigr)_{\infty}}\right\}
q^z\exp(-2\pi iz)
\end{split}
\end{equation*}
for some nonzero constant $C$. With similar arguments as in the proof of
Lemma \ref{boundRam} and Corollary \ref{boundRamcor}
(or Lemma \ref{boundAW} and Corollary \ref{corboundAW}) we conclude
that for some $C_1\in\mathbb{R}_{>0}$,
\[|I(z)|\leq C_1\,\epsilon^{|z|}
\]
for $z\in \mathbb{C}$ satisfying $\hbox{Re}(z)\geq 0$ and $\hbox{Im}(z)\leq 0$
as well as for $z\in\mathbb{C}$ satisfying $\hbox{Re}(z)\leq 0$ and 
$\hbox{Im}(z)\geq 0$,
with growth exponent $0<\epsilon=\exp(-2\pi m)<1$,
\[m=\min_{\theta\in [-\frac{\pi}{2},0]}\bigl((\tau-1)e^{i\theta}\bigr)>0.
\]
Thus we may apply Cauchy's theorem to rotate clockwise the integration 
cycle $\mathbb{R}$ 
in \eqref{NRfusedformula} to $-i\mathbb{R}$ about the origin without
altering its evaluation. Relaxing the parameter constraints leads to the
following result.
\begin{thm}[Hyperbolic Nassrallah-Rahman integral]
Let $\tau\in\mathbb{C}$ with $\hbox{Re}(\tau)<0$ and 
$\hbox{Im}(\tau)\geq 0$, and suppose that the parameters
$\tau_j\in\mathbb{C}$ \textup{(}$j=0,\ldots,4$\textup{)}
satisfy
\begin{equation}\label{pardomainNRb}
\hbox{Re}(\tau_j)<1 \quad (j=0,\ldots,4),\qquad
\hbox{Re}(a-\tau^{-1})>4,
\end{equation}
where $a=\tau_0+\tau_1+\tau_2+\tau_3+\tau_4$. Then
\begin{equation}\label{hyperbolicNR}
\int_{-i\infty}^{i\infty}
\frac{\lbrack 1\pm 2z,a-4\pm z;\tau\rbrack_{\infty}}
{\prod_{j=0}^4\lbrack \tau_j\pm z;\tau\rbrack_{\infty}}\,dz=
-\frac{2q^{-\frac{1}{24}}
\widetilde{q}\,{}^{\frac{1}{24}}}{\sqrt{-i\tau}}\,
\frac{\prod_{j=0}^4\lbrack a-\tau_j-3;\tau\rbrack_{\infty}}
{\prod_{0\leq k<m\leq 4}\lbrack \tau_k+\tau_m-1;\tau\rbrack_{\infty}}.
\end{equation}
\end{thm}
\begin{proof}
For $\tau\in \mathbb{H}\cap\mathbb{C}_-$ satisfying $\hbox{Re}(\tau^{-1})<-1$
we have already argued that \eqref{hyperbolicNR} is valid for 
parameters $\tau_j$ satisfying the constraints \eqref{pardomainNRa}.
By analytic continuation we conclude that \eqref{hyperbolicNR}
is valid under the milder parameter constraints \eqref{pardomainNRb},
using similar arguments as in the proof of Theorem \ref{hyperbolicR}.

The condition $\hbox{Re}(\tau^{-1})<-1$ can be removed by analytic
continuation using the asympotics (uniform for $\tau$ in compacta) of
the $\tau$-shifted factorial, see \eqref{as1} and \eqref{as2}.
The extension of \eqref{hyperbolicNR} to $\tau\in \mathbb{R}_{<0}$
then follows as in the proof of Theorem \ref{hyperbolicR}.
\end{proof}
\begin{rema}\label{limitNRAW}
The hyperbolic Askey-Wilson integral \eqref{AW3} is formally
a limit case of the hyperbolic Nassrallah-Rahman integral. For example,
let $\tau\in \mathbb{R}_{<0}$ and choose parameters $\tau_j\in \mathbb{C}$
($j=1,\ldots,4$) satisfying \eqref{pardomainAWb}. Concretely, 
the $\tau_j$'s thus satisfy
\[\hbox{Re}(\tau_j)<1,\qquad \hbox{Re}(\tau_1+\tau_2+\tau_3+\tau_4)>
3+\tau^{-1}.
\]
Choose $\tau_0^\prime\in\mathbb{R}_{<1}$ such that 
\[\tau_0^\prime+\hbox{Re}(\tau_1+\tau_2+\tau_3+\tau_4)>4+\tau^{-1},
\]
then any five tuple $(\tau_0,\tau_1,\tau_2,\tau_3,\tau_4)$ 
with $\hbox{Re}(\tau_0)=\tau_0^\prime$ satisfies the parameter constraints
\eqref{pardomainNRb}. We can now formally take the limit 
$\hbox{Im}(\tau_0)\rightarrow -\infty$ in the hyperbolic 
Nassrallah-Rahman integral \eqref{hyperbolicNR} while keeping
$\hbox{Re}(\tau_0)=\tau_0^\prime$ fixed. By \eqref{as3}, 
this formal limit of the hyperbolic Nassrallah-Rahman integral is the
hyperbolic Askey-Wilson integral \eqref{AW3}.
\end{rema}


\subsection{The hyperbolic degeneration of the elliptic
Nassrallah-Rahman integral}

We introduce the short hand notation $\underline{t}=(t_0,t_1,t_2,t_3,t_4)$.
With the notations of \S 1, we write
\begin{equation}
\Delta(z;\underline{t};p_1,p_2)=
\frac{\prod_{j=0}^4\Gamma(t_jz^{\pm 1};p_1,p_2)}
{\Gamma (z^{\pm 2}, Az^{\pm 1};p_1,p_2)}
\end{equation}
with $A=t_0t_1t_2t_3t_4$ and
\begin{equation}
N(\underline{t};p_1,p_2)=
\frac{2}{\bigl(p_1;p_1\bigr)_{\infty}\bigl(p_2;p_2\bigr)_{\infty}}
\frac{\prod_{0\leq k<m\leq 4}\Gamma(t_kt_m;p_1,p_2)}
{\prod_{j=0}^4\Gamma(At_j^{-1};p_1,p_2)},
\end{equation}
so that 
\begin{equation}\label{ell1}
\frac{1}{2\pi i}\int_{\mathbb{T}}
\Delta(z;\underline{t};p_1,p_2)\,\frac{dz}{z}=N(\underline{t};p_1,p_2)
\end{equation}
is Spiridonov's \cite{Sp}
elliptic Nassrallah-Rahman integral \eqref{betaEulerelliptic}, valid
for $|t_j|, |p_k|<1$ and $|p_1p_2|<|A|$.

Fix $\tau\in\mathbb{R}_{<0}$.
We introduce the limiting parameter $r\in\mathbb{R}_{>0}$ by writing
$\Delta_r(x;\underline{\tau};\tau)$ and 
$N_r(\underline{\tau};\tau)$
for the weight function $\Delta(\exp(2\pi irx);\underline{t};p_1,p_2)$
and the norm $N(\underline{t};p_1,p_2)$ with parameters
\[t_j=\exp(2\pi r(\tau_j-1)),\qquad
p_1=\exp(2\pi r/\tau),\qquad
p_2=\exp(-2\pi r).
\]
The conditions $|t_j|<1$ and $|p_1p_2|<|A|$ then translate to the
parameter constraints \eqref{pardomainNRb} for the corresponding 
hyperbolic Nassrallah-Rahman integral \eqref{hyperbolicNR}.

The elliptic beta integral
\eqref{ell1} can then be rewritten as
\begin{equation}\label{ell2}
\int_{-\frac{1}{2r}}^{\frac{1}{2r}}
\Delta_r(x;\underline{\tau};\tau)\,dx=
r^{-1}N_r(\underline{\tau};\tau).
\end{equation}
We show in this subsection that the formal limit $r\searrow 0$ of
\eqref{ell2} gives the hyperbolic Nassrallah-Rahman
integral \eqref{betaEuler5}.
The limit is based on Ruijsenaars' \cite[Prop. III.12]{R1}
observation that the hyperbolic gamma function is a limit case of the elliptic
gamma function (see \eqref{ellipticnormalized} and 
\eqref{ellhyplimit} for the explicit limit transition
in the present notations).

We first rewrite the weight function $\Delta_r(x;\underline{\tau};\tau)$
and the norm $N_r(\underline{\tau};\tau)$ in terms of the renormalized
elliptic gamma function $\widetilde{\Gamma}_r(z;\tau)$ (see
\eqref{ellipticnormalized}). This
yields the expressions
\begin{equation*}
\begin{split}
\Delta_r(x;\underline{\tau};\tau)&=
\frac{\prod_{j=0}^4\widetilde{\Gamma}_r\bigl(i-i\tau_j\pm x;\tau\bigr)}
{\widetilde{\Gamma}_r\bigl(\pm 2x,5i-ia\pm x;\tau\bigr)}
\,\exp\left(\frac{\pi}{2r}(1-\tau)\right),\\
N_r(\underline{\tau};\tau)&=
\frac{2}{\bigl(\exp(-2\pi r);\exp(-2\pi r)\bigr)_{\infty}
\bigl(\exp(2\pi r/\tau);\exp(2\pi r/\tau)\bigr)_{\infty}}\\
&\times\frac{\prod_{0\leq k<m\leq 4}
\widetilde{\Gamma}_r\bigl(2i-i\tau_k-i\tau_m;\tau\bigr)}
{\prod_{j=0}^4\widetilde{\Gamma}_r\bigl(4i-ia+i\tau_j;\tau\bigr)}\,
\exp\left(\frac{5\pi}{12r}(1-\tau)\right)
\end{split}
\end{equation*}
with $a=\tau_0+\tau_1+\tau_2+\tau_3+\tau_4$, where we use the notations
\begin{equation*}
\begin{split}
\widetilde{\Gamma}_r(a_1,\ldots,a_m;\tau)&=
\prod_{j=1}^m\widetilde{\Gamma}_r(a_j;\tau),\\
\widetilde{\Gamma}_r(a_1\pm z_1,\ldots,a_m\pm z_m;\tau)&=
\prod_{j=1}^m\widetilde{\Gamma}_r(a_j+z_j,a_j-z_j;\tau)
\end{split}
\end{equation*}
for products of renormalized elliptic gamma functions. 
Consequently,
\eqref{ell2} can be rewritten as
\begin{equation}\label{ell3}
\int_{-\frac{1}{2r}}^{\frac{1}{2r}}
\frac{\prod_{j=0}^4\widetilde{\Gamma}_r\bigl(i-i\tau_j\pm x;\tau\bigr)}
{\widetilde{\Gamma}_r\bigl(\pm 2x,5i-ia\pm x;\tau\bigr)}\,dx
=c_r(\tau)\,\frac{\prod_{0\leq k<m\leq 4}
\widetilde{\Gamma}_r\bigl(2i-i\tau_k-i\tau_m;\tau\bigr)}
{\prod_{j=0}^4\widetilde{\Gamma}_r\bigl(4i-ia+i\tau_j;\tau\bigr)}
\end{equation}
with 
\[c_r(\tau)=\frac{2\exp\left(\frac{\pi}{12r}(\tau-1)\right)}
{r\,\bigl(\exp(-2\pi r);\exp(-2\pi r)\bigr)_{\infty}
\bigl(\exp(2\pi r/\tau);\exp(2\pi r/\tau)\bigr)_{\infty}}.
\]
The limit $r\searrow 0$ of the constant $c_r(\tau)$ can
be computed using the modularity of the Dedekind eta function.
\begin{lem}\label{Ded}
\[\lim_{r\searrow 0}c_r(\tau)=\frac{2}{\sqrt{-\tau}}.
\]
\end{lem}
\begin{proof}
The Dedekind eta function $\eta(\sigma)$ for $\sigma\in\mathbb{H}$
is defined by
\begin{equation}\label{Dedekind}
\eta(\sigma)=\bigl(\exp(2\pi i\sigma);\exp(2\pi i\sigma)\bigr)_{\infty}
\exp(\pi i\sigma/12),
\end{equation}
and satisfies
\begin{equation}\label{Dedekindinversion}
\eta\left(-\frac{1}{\sigma}\right)=\eta(\sigma)\sqrt{-i\sigma}.
\end{equation}
Applying \eqref{Dedekindinversion} 
with $\sigma=ir$ and with $\sigma=-ir/\tau$,
we can rewrite $c_r(\tau)$ as
\[c_r(\tau)=\frac{2\exp\left(\frac{\pi r}{12\tau}(1-\tau)\right)}
{\sqrt{-\tau}\,\bigl(\exp(-2\pi/r);\exp(-2\pi/r)\bigr)_{\infty}
\bigl(\exp(2\pi\tau/r);\exp(2\pi\tau/r)\bigr)_{\infty}}.
\]
Since $\tau\in\mathbb{R}_{<0}$ we can take the limit $r\searrow 0$ in
the latter expression, which yields the desired result.
\end{proof}
We now further investigate the limit $r\searrow 0$ of 
\eqref{ell3} using \eqref{ellhyplimit}. We associate with $\tau$
the two deformation parameters $q=\exp(2\pi i\tau)$ and 
$\widetilde{q}=\exp(-2\pi i/\tau)$, and we use the standard
convention that $q^u=\exp(2\pi iu\tau)$ and 
$\widetilde{q}\,{}^u=\exp(-2\pi iu/\tau)$ for $u\in\mathbb{C}$.
For the weight function and norm 
in \eqref{ell3} we obtain by \eqref{ellhyplimit}
\begin{equation}\label{ad}
\begin{split}
\lim_{r\searrow 0}\,\frac{\prod_{j=0}^4
\widetilde{\Gamma}_r\bigl(i-i\tau_j\pm x;\tau\bigr)}
{\widetilde{\Gamma}_r\bigl(\pm 2x, 5i-ia\pm x;\tau\bigr)}
&=\alpha\,\frac{\lbrack 1\pm 2ix,a-4\pm ix;\tau\rbrack_{\infty}}
{\prod_{j=0}^4\lbrack \tau_j\pm ix;\tau\rbrack_{\infty}},\\
\lim_{r\searrow 0}\,\frac{\prod_{0\leq k<m\leq 4}
\widetilde{\Gamma}_r\bigl(2i-i\tau_k-i\tau_m;\tau\bigr)}
{\prod_{j=0}^4\widetilde{\Gamma}_r\bigl(4i-ia+i\tau_j;\tau\bigr)}
&=\beta\,\frac{\prod_{j=0}^4\lbrack a-\tau_j-3;\tau\rbrack_{\infty}}
{\prod_{0\leq k<m\leq 4}\lbrack \tau_k+\tau_m-1;\tau\rbrack_{\infty}}
\end{split}
\end{equation}
with the constants $\alpha$ and $\beta$ given by
\begin{equation*}
\begin{split}
\alpha&=q^{-\frac{1}{8}}\widetilde{q}\,{}^{\frac{1}{8}}
q^{-\frac{1}{2}\left\{\left(\frac{1}{2\tau}-\frac{1}{2}\right)^2
+\left(\frac{1}{2\tau}+\frac{9}{2}-a\right)^2
-\sum_{j=0}^4\left(\frac{1}{2\tau}+\frac{1}{2}-
\tau_j\right)^2\right\}},\\
\beta&=q^{-\frac{5}{48}}\widetilde{q}\,{}^{\frac{5}{48}}
q^{\frac{1}{4}\left\{\sum_{0\leq k<m\leq 4}
\left(\frac{1}{2\tau}+\frac{3}{2}-\tau_k-\tau_m\right)^2
-\sum_{j=0}^4\left(\frac{1}{2\tau}+\frac{7}{2}-a+\tau_j\right)^2
\right\}}.
\end{split}
\end{equation*}
Writing out the squares leads to 
\begin{equation}\label{ratio}
\beta=q^{-\frac{1}{24}}\widetilde{q}\,{}^{\frac{1}{24}}
\exp\left(\frac{\pi i}{4}\right)\alpha.
\end{equation}
Combining Lemma \ref{Ded}, \eqref{ad} and \eqref{ratio},
we can formally take the limit $r\searrow 0$ in \eqref{ell3},
leading to
\begin{equation}\label{ell4}
\int_{-\infty}^{\infty}\frac{\lbrack 1\pm 2ix,a-4\pm ix;\tau\rbrack_{\infty}}
{\prod_{j=0}^4\lbrack \tau_j\pm ix;\tau\rbrack_{\infty}}\,dx=
\frac{2q^{-\frac{1}{24}}\widetilde{q}\,{}^{\frac{1}{24}}}
{\exp(-\pi i/4)
\sqrt{-\tau}}\,\frac{\prod_{j=0}^4\lbrack a-\tau_j-3;\tau\rbrack_{\infty}}
{\prod_{0\leq k<m\leq 4}\lbrack \tau_k+\tau_m-1;\tau\rbrack_{\infty}}.
\end{equation}
Changing the integration variable and using
\[\exp(-3\pi i/4)\sqrt{-\tau}=-\sqrt{-i\tau}
\]
we arrive at
\[\int_{-i\infty}^{i\infty}\frac{\lbrack 1\pm 2x,a-4\pm x;\tau\rbrack_{\infty}}
{\prod_{j=0}^4\lbrack \tau_j\pm x;\tau\rbrack_{\infty}}\,dx=
-\frac{2q^{-\frac{1}{24}}\widetilde{q}\,{}^{\frac{1}{24}}}
{\sqrt{-i\tau}}\,\frac{\prod_{j=0}^4\lbrack a-\tau_j-3;\tau\rbrack_{\infty}}
{\prod_{0\leq k<m\leq 4}\lbrack \tau_k+\tau_m-1;\tau\rbrack_{\infty}},
\]
which is the hyperbolic Nassrallah-Rahman integral \eqref{betaEuler5}.

\section{Appendix: The hyperbolic gamma function}

In this section we discuss the $\tau$-shifted factorial 
$\lbrack z;\tau\rbrack_\infty$ (see \eqref{taufactorial})
and its connection with
Ruijsenaars' \cite{R1} hyperbolic gamma function.
For detailed proofs we refer to Ruijsenaars' papers 
\cite{R1} and \cite{R2}.

We start by recalling the definition of Ruijsenaars' hyperbolic gamma 
function. We write 
\[\mathbb{C}_{\pm}=\{z\in \mathbb{C} \, | \, \hbox{Re}(z)\gtrless 0 \}
\] 
for the open right/left half plane.
The integral
\begin{equation}\label{g}
g(z)=g(a_+,a_-;z)=
\int_0^{\infty}\frac{dy}{y}
\left(\frac{\sin(2yz)}{2\sinh(a_+y)\sinh(a_-y)}-\frac{z}{a_+a_-y}\right)
\end{equation}
defines an analytic function for $(a_+,a_-,z)\in\mathcal{D}$, where
\[\mathcal{D}=\{(a_+,a_-,z)\in\mathbb{C}^{\times 3} \, | \,
\hbox{Re}(a_{\pm})>0,\,\,\, |\hbox{Im}(z)|<\frac{1}{2}(\hbox{Re}(a_+)+
\hbox{Re}(a_-)) \}.
\]
Ruijsenaars' \cite{R1}
{\it hyperbolic gamma function} is now defined 
as the analytic, zero-free function on $\mathcal{D}$ defined by 
\begin{equation}\label{G}
\Gamma_h(z)=\Gamma_h(a_+,a_-;z)=\exp(ig(a_+,a_-;z)).
\end{equation}
The hyperbolic gamma function satisfies the functional equations
\begin{equation}\label{ADEh}
\frac{\Gamma_h(z+ia_\pm/2)}{\Gamma_h(z-ia_\pm/2)}=2\cosh(\pi z/a_\mp)
\end{equation}
whenever the left hand side is defined. Hence 
$\Gamma_h$ extends to a meromorphic function
in the domain $(a_+,a_-,z)\in \mathbb{C}_+^{\times 2}\times 
\mathbb{C}$, which we again denote by $\Gamma_h(a_+,a_-;z)$.

The functional equations \eqref{ADEh} imply that the zeros of 
$z\mapsto \Gamma_h(a_+,a_-;z)$ are located at
\[\left(\mathbb{Z}_{\geq 0}+\frac{1}{2}\right)ia_++
\left(\mathbb{Z}_{\geq 0}+\frac{1}{2}\right)ia_-,\]
and the poles are located at 
\[-\left(\mathbb{Z}_{\geq 0}+\frac{1}{2}\right)ia_+
-\left(\mathbb{Z}_{\geq 0}+\frac{1}{2}\right)ia_-.\]
The zeros and poles are simple when 
$a_+/a_-\not\in\mathbb{Q}_{>0}$. 

The explicit integral expression for $g$ implies 
\begin{equation}\label{aut1}
\Gamma_h(a_+,a_-;z)=\Gamma_h(a_-,a_+;z),
\end{equation}
\begin{equation}\label{aut2}
\Gamma_h(ra_+,ra_-;rz)=\Gamma_h(a_+,a_-;z), \qquad r\in\mathbb{R}_{>0},
\end{equation}
as well as the {\it reflection equation}
\begin{equation}\label{aut3}
\Gamma_h(a_+,a_-;z)\Gamma_h(a_+,a_-;-z)=1.
\end{equation}
Note that by the translation invariance
\eqref{aut2}, the hyperbolic gamma function $\Gamma_h$ essentially only
depends on the quotient $a_+/a_-$ of the two deformation parameters $a_{\pm}$.

Some special values of $\Gamma_h$ are easily computed. Obviously, we 
have 
\begin{equation}\label{gamma0}
\Gamma_h(a_+,a_-;0)=1.
\end{equation}
Applying the two functional equations and the reflection equation, 
we furthermore have
\[\Gamma_h\left(a_+,a_-;z-\frac{ia_+}{2}+\frac{ia_-}{2}\right)
\Gamma_h\left(a_+,a_-;-z-\frac{ia_+}{2}+\frac{ia_-}{2}\right)=
\frac{\sinh(\pi z/a_+)}{\sinh(\pi z/a_-)}.
\]
Taking the limit $z\rightarrow 0$ we
conclude that
\begin{equation}\label{gamma1}
\Gamma_h\left(a_+,a_-;\frac{ia_-}{2}-\frac{ia_+}{2}\right)=
\sqrt{\frac{a_-}{a_+}}
\end{equation}
(to see that the branch of the square root is the right one, note that
$\Gamma_h(a_+,a_-;x)>0$ for $a_\pm\in\mathbb{R}_{>0}$ 
and $x\in i\mathbb{R}$ in 
view of \eqref{g} and \eqref{G}).

The hyperbolic gamma function can be explicitly expressed as quotient
of trigonometric gamma functions
when $\hbox{Im}(a_+/a_-)\not=0$. 
This expression was first
obtained by Shintani \cite{Sh}, see also \cite[Appendix A]{R2}.
\begin{prop}\label{Shintani}
Let $a_{\pm}\in\mathbb{C}_+$ with $\hbox{Im}(a_+/a_-)>0$.
Then
\begin{equation*}
\begin{split}
\Gamma_h(a_+,a_-;z)=&\frac{\bigl(-\exp(\pi ia_+/a_-)\exp(-2\pi z/a_-);
\exp(2\pi ia_+/a_-)\bigr)_{\infty}}
{\bigl(-\exp(-\pi ia_-/a_+)\exp(-2\pi z/a_+);
\exp(-2\pi ia_-/a_+)\bigr)_{\infty}}\\
&\times\exp\left(-\frac{\pi i}{24}
\left(\frac{a_+}{a_-}+\frac{a_-}{a_+}\right)\right)\,
\exp\left(-\frac{\pi iz^2}{2a_+a_-}\right).
\end{split}
\end{equation*}
\end{prop}
\begin{proof}
We sketch a proof because the infinite product expression
for $\Gamma_h$ plays such a crucial role in the present paper.
Write $\widehat{\Gamma}_h(z)$
for the right hand side of the desired identity. It is easily verified
that $\widehat{\Gamma}_h(z)$ is meromorphic in $z$, having the same
poles and zeros as $\Gamma_h(z)$. 
A direct check shows that $\widehat{\Gamma}_h(z)$
satisfies the same functional equations \eqref{ADEh}
as $\Gamma_h(z)$. Thus $\widehat{\Gamma}_h/\Gamma_h$ is an entire, bounded
function, hence a constant. The constant is one since
\[\widehat{\Gamma}_h\left(\frac{ia_-}{2}-\frac{ia_+}{2}\right)=
\sqrt{\frac{a_-}{a_+}},
\]
which follows from the modularity \eqref{Dedekindinversion}
of the Dedekind eta function \eqref{Dedekind}.
\end{proof}
Note that an infinite product expression for $\Gamma_h(a_+,a_-;z)$ 
with $\hbox{Im}(a_+/a_-)<0$ can be directly derived from Proposition
\ref{Shintani} by applying \eqref{aut1}.

Using Proposition \ref{Shintani}, 
the rather harmless looking reflection equation 
\eqref{aut3} turns into a 
nontrivial infinite product identity.
This identity can be proven without referring to the integral
representation of $\Gamma_h$ using Jacobi's 
triple product identity, Jacobi's inversion formula and the modularity
of the Dedekind eta function.

Ruijsenaars \cite{R2} established rather delicate asymptotic bounds for the
hyperbolic gamma function. We formulate here a weaker version of 
these bounds which suffices for our purposes. It can be stated as 
\[\Gamma_h(a_+,a_-;z)=\mathcal{O}\left(\exp\left(\mp 
\frac{i\pi z^2}{2a_+a_-}\right)\right),
\qquad \hbox{Re}(z)\rightarrow \pm \infty,
\]
uniformly for $\hbox{Im}(z)$ in compacta of $\mathbb{R}$ and 
for $a_\pm$ in compacta of $\mathbb{C}_+$. The precise meaning is as follows,
cf. \cite[Thm. A.1]{R2}.
\begin{prop}\label{asympt1}
Let $K_{\pm}\subset \mathbb{C}_+$ and $K\subset \mathbb{R}$ be compact subsets.
There exist positive constants $R=R(K_+,K_-;K)$ and 
$C=C(K_+,K_-;K)$, both depending on $K_{\pm}$ and $K$ only,
such that
\[\left|\Gamma_h(a_+,a_-;z)\exp\left(\pm \frac{i\pi z^2}{2a_+a_-}\right)\right|
\leq C,\qquad \hbox{Re}(z)\gtrless R
\]
when $\hbox{Im}(z)\in K$ and $a_{\pm}\in K_{\pm}$.
\end{prop}
The precise connection between $\Gamma_h$ and the $\tau$-shifted factorial
$\lbrack z;\tau\rbrack_{\infty}$ (see \eqref{taufactorial}) 
is as follows.
Fix $\tau\in \mathbb{C}_-\cap\mathbb{H}$
and write $q=\exp(2\pi i\tau)$ and $\widetilde{q}=\exp(-2\pi i/\tau)$.
Recall the notational convention $q^u=\exp(2\pi i\tau u)$ and 
$\widetilde{q}\,{}^u=\exp(-2\pi iu/\tau)$ for $u\in\mathbb{C}$.
The {\it $\tau$-shifted factorial} $\lbrack z;\tau\rbrack_{\infty}$
(see \eqref{taufactorial}) can then be expressed in terms
of the hyperbolic gamma function by
\begin{equation}\label{tauhyper}
\lbrack z;\tau\rbrack_{\infty}=
q^{\bigl(z-\frac{1}{2}-\frac{1}{2\tau}\bigr)^2/4}q^{-\frac{1}{48}}
\widetilde{q}\,{}^{\frac{1}{48}}\Gamma_h\left(-\frac{1}{\tau},1;
i\left(z-\frac{1}{2}-\frac{1}{2\tau}\right)\right).
\end{equation}
Indeed, substituting the infinite product expression for $\Gamma_h$
from Proposition \ref{Shintani} in the right hand side of \eqref{tauhyper}
gives precisely the expression \eqref{taufactorial} defining the
$\tau$-shifted factorial. Working with
the $\tau$-factorial $\lbrack z;\tau\rbrack_{\infty}$ has the advantage
that formulas have a similar appearance as in the 
usual basic hypergeometric (trigonometric) setup.
The disadvantage
is the loss of symmetry in the parameters $a_{\pm}$ (see \eqref{aut1}).
We end this section by reformulating the above properties of
the hyperbolic gamma function in terms of the $\tau$-shifted factorial.

Formula \eqref{tauhyper} shows that 
$\lbrack z;\tau\rbrack_{\infty}$ has a meromorphic continuation
to $(z,\tau)\in \mathbb{C}\times \mathbb{C}_-$,
which we again denote by $\lbrack z;\tau\rbrack_{\infty}$
(in fact, in view of \eqref{taufactorial} and \eqref{aut1}, 
$\lbrack z;\tau\rbrack_{\infty}$ extends to a meromorphic function 
in $(z,\tau)\in
\mathbb{C}\times\bigl(\mathbb{C}\setminus \mathbb{R}_{\geq 0}\bigr)$). 
The functional equations \eqref{ADEh} become
\begin{equation}\label{functionalequations}
\begin{split}
\lbrack z+1;\tau\rbrack_{\infty}&=\bigl(1-q^z\bigr)
\lbrack z;\tau\rbrack_{\infty},\\
\lbrack z-\tau^{-1};\tau\rbrack_{\infty}&=
\bigl(1-\widetilde{q}\,{}^{-1}\exp(-2\pi iz)\bigr)
\lbrack z;\tau\rbrack_{\infty}.
\end{split}
\end{equation}
The zeros of $z\mapsto\lbrack z;\tau\rbrack_{\infty}$ are located at
\begin{equation}\label{zero}
\frac{1}{\tau}\,\mathbb{Z}_{\leq 0}+\mathbb{Z}_{>0},
\end{equation}
and the poles are located at
\begin{equation}\label{pole}
\frac{1}{\tau}\,\mathbb{Z}_{>0}+\mathbb{Z}_{\leq 0}.
\end{equation}
The zeros and poles are simple when $\tau\not\in \mathbb{Q}_{<0}$.
Property \eqref{aut1} translates into
\begin{equation}
\lbrack z;\tau^{-1}\rbrack_{\infty}=
\lbrack -\frac{z}{\tau}+\frac{1}{\tau}+1;\tau\rbrack_{\infty}
\end{equation}
for $\tau\in \mathbb{C}_-$. The reflection equation \eqref{aut3}
becomes
\begin{equation}\label{reflectioneq}
\lbrack \frac{1}{2}+\frac{1}{2\tau}+z,
\frac{1}{2}+\frac{1}{2\tau}-z;\tau\rbrack_{\infty}=
q^{-\frac{1}{24}}\widetilde{q}\,{}^{\frac{1}{24}}q^{\frac{z^2}{2}}
\end{equation}
for $\tau\in\mathbb{C}_-$. Finally, the asymptotic bounds for 
the $\tau$-shifted factorial, deduced from Proposition \ref{asympt1},
become
\begin{equation}\label{as1}
\lbrack z;\tau\rbrack_{\infty}=
\mathcal{O}\left(q^{\bigl(z-\frac{1}{2}-\frac{1}{2\tau}\bigr)^2/2}\right),
\qquad \hbox{Im}(z)\rightarrow \infty,
\end{equation}
uniformly for $\hbox{Re}(z)$ in compacta of $\mathbb{R}$ and for
$\tau$ in compacta of $\mathbb{C}_-$, and
\begin{equation}\label{as2}
\lbrack z;\tau\rbrack_{\infty}=\mathcal{O}(1),
\qquad \hbox{Im}(z)\rightarrow -\infty,
\end{equation}
uniformly for $\hbox{Re}(z)$ in compacta of $\mathbb{R}$ and for
$\tau$ in compacta of $\mathbb{C}_-$. Using the more
precise asymptotic estimates for the hyperbolic gamma function in 
\cite[Thm. A.1]{R1}, we actually have the limit
\begin{equation}\label{as3}
\lim_{\hbox{Im}(z)\rightarrow -\infty}\lbrack z;\tau\rbrack_{\infty}=1
\end{equation}
uniformly for $\hbox{Re}(z)$ in compacta of $\mathbb{R}$ 
and for $\tau$ in compacta of $\mathbb{C}_-$.

We end this section by rewriting the hyperbolic degeneration 
(see \cite[Prop. III.12]{R1}) of Ruijsenaars' elliptic
gamma function in our present
notations. Ruijsenaars' elliptic gamma function
$G(r,a_+,a_-;z)$ (see e.g. \cite[Prop. III.11]{R1})
relates to the elliptic gamma function \eqref{ellgamma}
by
\[G(r,a_+,a_-;z)=\Gamma\bigl(\exp(2irz-a_+r-a_-r);\exp(-2a_+r),
\exp(-2a_-r)\bigr)
\]
for $\hbox{Re}(a_\pm r)>0$. Thus \cite[Prop. III.12]{R1} becomes
\begin{equation*}
\begin{split}
\lim_{r\searrow 0}\Gamma\bigl(\exp(2\pi irz);\exp(-2\pi a_+r),
\exp(-2\pi a_-r)\bigr)&\exp\left(\frac{2\pi z-\pi ia_+-\pi ia_-}
{12 ira_+a_-}\right)\\
&=\Gamma_h\left(a_+,a_-;z-\frac{ia_+}{2}-\frac{ia_-}{2}\right)
\end{split}
\end{equation*}
for $a_{\pm}\in\mathbb{R}_{>0}$.
We now take $a_+=-\frac{1}{\tau}$ and $a_-=1$ with
$\tau\in \mathbb{R}_{<0}$, and we denote for simplicity
\begin{equation}\label{ellipticnormalized}
\widetilde{\Gamma}_r\bigl(z;\tau\bigr)=
\Gamma\bigl(\exp(2\pi irz);\exp(2\pi r/\tau),\exp(-2\pi r)\bigr) 
\exp\left(\frac{\pi i(2\tau z+i-i\tau)}{12r}\right).
\end{equation}
By \eqref{tauhyper}, we arrive at the limit
\[
\lim_{r\searrow 0}\widetilde{\Gamma}_r(z;\tau)=
q^{-\bigl(\frac{1}{2\tau}-\frac{1}{2}-iz\bigr)^2/4}
q^{\frac{1}{48}}\widetilde{q}\,{}^{-\frac{1}{48}}\lbrack \tau^{-1}-iz;
\tau\rbrack_{\infty}
\]
with $q=\exp(2\pi i\tau)$ and $\widetilde{q}=\exp(-2\pi i/\tau)$.
Applying the reflection equation \eqref{reflectioneq} leads to the limit
\begin{equation}\label{ellhyplimit}
\lim_{r\searrow 0}\widetilde{\Gamma}_r(z;\tau)=
q^{\bigl(\frac{1}{2\tau}-\frac{1}{2}-iz\bigr)^2/4}
q^{-\frac{1}{48}}\widetilde{q}\,{}^{\frac{1}{48}}\lbrack 1+iz;
\tau\rbrack_{\infty}^{-1}
\end{equation}
for $\tau\in \mathbb{R}_{<0}$.


\end{document}